\documentclass[11pt]{article}
\usepackage[short]{optidef}
\usepackage{tikz, lscape, amsmath,amsfonts,latexsym,amssymb}
\usepackage{amsthm}
\usepackage{thmtools}

\usepackage{pgf} \usepackage{tikz}
\usetikzlibrary{calc} \usetikzlibrary{arrows}
\usetikzlibrary{shapes} %
\usetikzlibrary{shapes}
\usetikzlibrary{plotmarks}
\usepackage{accents}

\usepackage[utf8]{inputenc}
\usepackage[english]{babel}
\usepackage{lmodern, verbatim}
\usepackage[T1]{fontenc}
\usepackage{calc,bm}
\usepackage{tabularx,enumitem,nameref,hyperref}
\usepackage[mathscr]{eucal}
\usepackage{subfig,stackrel,floatrow,xspace}
\usepackage{mathtools}

\usepackage{cleveref}

\usepackage[margin=1in]{geometry}%

\newcommand{\NN}{\mathbb{N}}

\newcommand{\RR}{\mathbb{R}}

\newcommand{\mR}{\mathbb{R}_\bot}
\newcommand{\pR}{\mathbb{R}_\top}
\newcommand{\pmR}{\overline{\mathbb{R}}}
\newcommand{\pmRX}{\pmR^\X}
\newcommand{\pmRY}{\pmR^\Y}

\newcommand{\pmRmax}{\pmR_{\max}}

\newcommand{\G}{\mathcal{G}}

\newcommand{\bB}{\bar B}
\newcommand{\bF}{\bar F}
\newcommand{\xb}{\bar x}

\newcommand{\bp}{\bar p}

\newcommand{\bJ}{\bar J}
\newcommand{\bGinf}{\overline \G^{\inf}}

\newcommand{\rfunk}{\operatorname{RFunk}}

\newcommand{\hilbert}{\operatorname{Hil}}
\newcommand{\interior}{\operatorname{int}}

\newcommand{\Asc}{\mathscr A}
\newcommand{\Lcal}{\mathcal{L}}

\newcommand{\infG}{\inf\nolimits^\G}
\newcommand{\infnl}{\inf\nolimits}

\makeatletter
\newcommand{\dotminus}{\mathbin{\text{\@dotminus}}}

\newcommand{\@dotminus}{%
	\ooalign{\hidewidth\raise1ex\hbox{.}\hidewidth\cr$\m@th-$\cr}%
}
\newcommand{\dotbminus}{\mathbin{\text{\@dotbminus}}}

\newcommand{\@dotbminus}{%
	\ooalign{\hidewidth\raise-0.2ex\hbox{.}\hidewidth\cr$\m@th-$\cr}%
}

\newcommand{\dotbplus}{\mathbin{\text{\@dotbplus}}}

\newcommand{\@dotbplus}{%
	\ooalign{\hidewidth\raise-0.6ex\hbox{.}\hidewidth\cr$\m@th+$\cr}%
}

\makeatother

\newcommand{\A}{\mathcal{A}}

\newcommand{\C}{\mathcal{C}}

\newcommand{\X}{\mathscr{X}}
\newcommand{\Y}{\mathscr{Y}}
\newcommand{\Z}{\mathscr{Z}}

\newcommand{\la}{\lambda}
\newcommand{\Min}{\operatorname{Min}}

\newcommand{\Gsup}{\overline{\G}^{\sup}}
\newcommand{\Fsup}{\overline{\F}^{\sup}}

\DeclareMathOperator*{\Sp}{span}
\DeclareMathOperator*{\Rg}{Rg}

\DeclareMathOperator*{\Lipb}{\operatorname{Lip}_1}

\DeclareMathOperator*{\Dom}{Dom}

\DeclareMathOperator*{\Id}{Id}
\newcommand{\F}{\mathscr{F}} %

\newcommand{\tb}{\textbf} %
\newcommand{\R}{\RR} %
\newcommand{\Q}{\mathbb Q} %
\newcommand{\N}{\NN} %

\newcommand{\bracket}[1]{\langle #1 \rangle}
\newcommand{\bk}[1]{\langle #1 \rangle}

\theoremstyle{definition}

\theoremstyle{plain}
\newtheorem{Theorem}{Theorem}
\newtheorem*{Theorem*}{Theorem}
\newtheorem{Assumption}{Assumption}
\newtheorem{Corollary}[Theorem]{Corollary}%
\newtheorem{Proposition}[Theorem]{Proposition}%
\newtheorem{Claim}[Theorem]{Claim}%
\newtheorem{Lemma}[Theorem]{Lemma}%

\theoremstyle{definition}
\newtheorem{Definition}[Theorem]{Definition}%
\newtheorem{Example}[Theorem]{Example}%
\newtheorem{Remark}[Theorem]{Remark}%
\date{\today}

\makeatletter
\let\orgdescriptionlabel\descriptionlabel
\renewcommand*{\descriptionlabel}[1]{%
	\let\orglabel\label
	\let\label\@gobble
	\phantomsection
	\edef\@currentlabel{#1\unskip}%
	\let\label\orglabel
	\orgdescriptionlabel{#1}%
}
\makeatother

\title{Order isomorphisms of sup-stable function spaces:\\ continuous, Lipschitz, c-convex, and beyond}
\author{Pierre-Cyril Aubin-Frankowski\footnote{CERMICS, ENPC, Institut Polytechnique de Paris, Marne-la-Vallée, France, pierre-cyril.aubin@enpc.fr (corresponding author)}, Stéphane Gaubert\footnote{Inria Saclay
CMAP, École polytechnique, IP Paris
France, stephane.gaubert@inria.fr}}%
\date{\today}%

\begin{document}
	\maketitle
 
    \begin{abstract}
      There have been many parallel streams of research studying order isomorphisms of some specific sets $\G$ of functions from a set $\X$ to $\R\cup\{\pm\infty\}$, such as the sets of convex or Lipschitz functions. We develop in this article a unified approach inspired by $c$-convex functions. Our results are obtained highlighting the role of inf and sup-irreducible elements of $\G$ and the usefulness of characterizing them, to subsequently derive the structure of order isomorphisms, and in particular of those commuting with the addition of scalars. We show that in many cases all these isomorphisms $J:\G\to\G$ are of the form $Jf=g+f\circ \phi$ for a translation $g:\X\to\R$ and a bijective reparametrization $\phi:\X\to \X$. Given a reference anti-isomorphism, this characterization then allows to recover all the other anti-isomorphisms. We apply our theory to the sets of $c$-convex functions on compact Hausdorff spaces, to the set of lower semicontinuous (convex) functions on a Hausdorff topological vector space and to 1-Lipschitz functions of complete metric spaces. The latter application is obtained using properties of the horoboundary of a metric space.\\


      \noindent MSC classes:	15A80, 06D50, 26B25, 06A15.
    \end{abstract}

    \paragraph{\bf Context.}  We are interested in the characterization of the order preserving (anti)-isomorphisms for max or (max,+) operations over some sets $\G$ of functions from a set $\X$ to $\R\cup\{\pm \infty\}$ such that $\G$ is stable by arbitrary suprema, and possibly by the addition of scalars. Why should one be interested in such characterization of transformations? First because it corresponds to nonlinear versions of the Banach-Stone theorem, and thus relates to studying morphisms between sets, as pursued in \cite{weaver1994lattices,Leung2016}. Second, the characterization shines a new light on the uniqueness of some very familiar transformations such as the Fenchel transform \cite{ArtsteinAvidan2009} or the polar duality of convex bodies \cite{ArtsteinAvidan2008,Brczky2008}, especially in the context of valuation theory in convex geometry \cite{Ludwig2023}. Finally, dualities, i.e.\ order reversing involutions, are symmetries for some equations like the complex Monge-Ampère \cite{Lempert2017,Berndtsson2020} and allow for new Blaschke-Santaló-type inequalities as shown by \cite{ArtsteinAvidan2023}. The latter two topics are related to the flourishing field of optimal transport, being also targeted by our study with its focus on the set of $c$-convex functions.

    \paragraph{\bf Main results.} We show that, for several function spaces $\G\subset \pmRX$ of interest, every bijection $J$ over $\G$ that commutes with the pairwise supremum and the addition of scalars, in short a (max,+)-isomorphism, is of the form \eqref{eq:charac_J_intro} owing to the following characterization:%
    \begin{Theorem}[Main theorem on (max,+)-isomorphisms]\label{thm:max_plus_isophi_general_intro}%
        Let $\X$ and $\Y$ be nonempty sets. Let $\G$ (resp.\ $\F$) be a subset of $\pmRX$ (resp.\ $\pmRY$), with both $\F$ and $\G$ being \emph{proper and separating}, stable by arbitrary suprema and addition of scalars. Let $J$ be a (max,+)-isomorphism from $\F$ onto $\G$ and define $e_x=\sup_{u\in\G}u(\cdot)-u(x)$ and $e'_y=\sup_{v\in\F}v(\cdot)-v(y)$. Then the following statements are equivalent:
		\begin{enumerate}[label=\roman*),labelindent=0cm,leftmargin=*,topsep=0.1cm,partopsep=0cm,parsep=0.1cm,itemsep=0.1cm]
			\item \label{it_J_e_x_intro} for all $y \in \Y$, $J(e'_y)\in\{e_x+\la\}_{x\in\X,\,\la\in\pmR}$;
			\item \label{it_g_f_phi_intro} there exists $g:\X\rightarrow \R$ and a bijective $\phi:\X \rightarrow \Y$ such that
            \begin{equation}\label{eq:charac_J_intro}
			Jf(x)=g(x)+f(\phi(x))
		\end{equation}
            and for all $f\in\F$, $h\in\G$, $(g+f\circ\phi)\in\G$ and $(-g\circ\phi^{-1}+h\circ\phi^{-1})\in\F$.
		\end{enumerate}
     Conversely every $J$ satisfying \ref{it_g_f_phi_intro} is a (max,+)-isomorphism, and, the representation $(g,\phi)$ in \eqref{eq:charac_J_intro} is unique. 
	\end{Theorem}
   The precise statement is given in \Cref{thm:max_plus_isophi_general}, \emph{proper and separating} corresponding to \Cref{assumption:proper}, and the functions $e_x$ being introduced in \eqref{eq:def_e_x}. The other properties of $g$ and $\phi$ depend obviously on the specific $\F$ and $\G$ studied, some characterizations obtained are collected in \Cref{tab:spaces_intro}. The particularity of our approach is that we encompass all these functions spaces within a single abstract formalism of sup-stable spaces. \Cref{thm:max_plus_isophi_general} establishing \eqref{eq:charac_J_intro} is only a posteriori adapted to characterize the $g$ and $\phi$ for each specific space.

    \begin{table}[!h]
		\caption{Some (max,+)-isomorphisms based on \eqref{eq:charac_J_intro} for specific spaces}\scriptsize
		\label{tab:spaces_intro}
		\begin{tabular}{@{}c|c|c|c|c@{}}
			Set $\X$ & Function space $\G$ & Translation $g$ & Reparametrization $\phi$ & \\ \hline
			Hausdorff topological space & \begin{tabular}{c}lower semicontinuous\\  functions\end{tabular}  & continuous & homeomorphism & \Cref{thm:max_plus_isophi_lsc} \\ \hline
			 \begin{tabular}{c}complete\\ metric space\end{tabular}  & 1-Lipschitz functions  & constant & isometry & \Cref{th-Lip}\\ \hline

			 \begin{tabular}{c}locally
        convex Hausdorff \\  topological vector space\end{tabular} & \begin{tabular}{c}lower semicontinuous\\ convex functions\end{tabular}  &
			continuous affine & \begin{tabular}{c}invertible affine\\with w*-w* adjoint\end{tabular} & \Cref{thm:isophi_convex}  
		\end{tabular}
	\end{table}
    Concerning order reversing isomorphisms $T$ between sets $\F$ and $\G$ commuting with the addition of scalars, we characterize their existence through kernels $b$ in \Cref{thm:sym_anti-involution}. Owing to the immediate \Cref{prop:isophi_anti-invol}, the uniqueness of $T$ up to basic transformations is instead a consequence of results like \Cref{thm:max_plus_isophi_general_intro}.
 
     On $c$-convex functions over a compact set, we show in \Cref{thm:rg_B_isometry} that the existence of a (max,+)-isomorphism between $b$-convex and $c$-convex functions is equivalent to a ``conjugacy'' formula \eqref{eq:kernel_isometry} satisfied by the kernels. Concerning 1-Lipschitz functions, the results rely on the compactification of metric spaces by horofunctions, and in particular on properties of Busemann functions (\Cref{th-agw}). We also allow for weak distances, i.e.\ potentially assymetric and with negative values.%
     
     Similarly to the proof of the Banach-Stone theorem, see e.g.\ the lecture notes \cite{Garrido2002}, the key idea in our proofs is to advocate for inf and sup-irreducibility as a good notion of ``extremality'' (\Cref{def:irreducible_points}) and to show that these ``extremal'' points are preserved by order isomorphisms (\Cref{lem:isophi_extreme_pts}) and generate the whole space, in a Krein--Milman fashion (\Cref{lem:e_x_and_Ginf} and \Cref{fct:extreme_pts_generators_minSx}).%

    Beyond (max,+)-isomorphisms, we study the order isomorphisms of sets containing the indicator functions of points, such as the set of lower semicontinuous functions on Hausdorff topological spaces (\Cref{thm:max_isophi_lsc}) and its subset of convex functions (\Cref{thm:isophi_convex}). In these cases, indicator functions play a major role in the characterization. We in particular obtain the extension of the theorems of \cite{ArtsteinAvidan2009,Iusem2015} in the generality of locally convex Hausdorff spaces.
    
    \paragraph{\bf Related work.} Some of our results have precursors when $\X$ is finite, in the setting of tropical geometry -- the (max,+) structure considered here being the main example of tropical algebraic structure. In particular, the notion of irreducibility was used in \cite{BSS,GK} when studying the max-plus analogue of the representation of convex sets in terms of extreme points or extreme rays.
    Concerning infinite sets, the literature can be grouped according to the function spaces studied:

    First, on convex functions, our focus meets that of the research program set by S.\ Artstein-Avidan and V.\ Milman on the characterization of some operations on subsets of convex functions. In their seminal paper \cite{ArtsteinAvidan2009}, they have shown that, for $\G$ the space of convex lower semicontinuous functions over $\R^d$, the only possible anti-involutions are compositions of the Fenchel transform with affine transformations. This result was further extended to Banach spaces by \cite{Iusem2015}, followed by \cite{Cheng2021} who characterized for which Banach spaces the Fenchel transform is an anti-involution. Using techniques specific to the convex case, a further generalization to convex subsets of locally convex Hausdorff spaces was obtained in \cite[Theorem 2.13]{leung2023}, which contains our \Cref{thm:isophi_convex}. We refer to \cite{Iusem2015,ArtsteinAvidan2023} for a list of other developments on $s$-concave functions or convex bodies. Underlying all of these works on convex functions is a fundamental result of affine geometry that, in dimensions larger than 2, transformations preserving straight lines are affine. The notable difference of our approach is that we consider spaces beyond convex functions on vector spaces, so in general we do not have linear structures.%
    
    Second, we tackle spaces of (lower semi) continuous functions. This is related to developing analogues of the Banach-Stone theorem for order preserving bijections, as for instance done for the set of continuous functions $C(\X,\R)$ on a compact set $\X$ in \cite{Kaplansky1948,Sanchez2008} or in the max-plus literature \cite[Chapter 2]{KM}. The isomorphisms of the set $C(\X,\R)$ and its subsets have indeed a long history based on many viewpoints, as nicely recalled in \cite[p2]{Leung2016}. Some interpretations can be more insightful, for instance, seeing $C(\X,\R)$ as a vector lattice, Kaplansky’s Theorem gives the famed Gelfand-Kolmogorov theorem which focuses on $C(\X,\R)$ as an algebra. Here we switch instead the focus to the sup-closure of $C(\X,\R)$, namely the space of lower semicontinuous functions. The isomorphisms of the latter have received less interest so far, and we highlight the role of the indicator functions.

    Finally, an important application of ours is to the set of Lipschitz functions on a complete metric space with Lipschitz constant at most 1. Most of the interest for spaces of Lipschitz functions has been for their Banach structure, but they also have a rich lattice structure as advocated by \cite{weaver1994lattices}. This induces quite some polysemy on the notion of isomorphism. For instance \cite{Candido2019} studies Banach space isomorphisms between spaces of Lipschitz functions. Weaver leveraged instead the inf/sup-lattice properties of the vector spaces of Lipschitz functions to characterize linear Banach space isomorphisms preserving Riesz-norms \cite{weaver1994lattices}. Banach–Stone-like theorems for Lipschitz functions have also been shown for linear isomorphisms such as \cite[Theorem 4.1]{Jimenez-Vargas2009} and nonlinear \cite[Theorem 1]{Sanchez2011} and references within. Finally Leung et al.\ for linear, in \cite{Leung2013}, and nonlinear order isomorphisms, in \cite{Leung2016}, designed an abstract framework to describe order isomorphisms on ``near vector lattices''. However all the aforementioned authors rely extensively on vector subspaces of the space of continuous functions, while we consider sup-stable sets of functions. The closest research to our work was done by \cite{Cabello2017,Daniilidis2020} for $1$-Lipschitz functions on metric spaces with non-symmetric distance with the assumption that the isomorphisms also preserve convex combinations. Nevertheless, none of the authors uses the notions of inf or sup-irreducible functions. Moreover, a central
    idea in our proof is the use of the compactification of a metric space by horofunctions. We exploit the fact
    that Busemann points (limits of geodesics) are precisely the non-trivial
    sup-irreducible elements of the space of $1$-Lipschitz functions~\cite{AGW}.
    
    There is a separate direction of research to characterize cones for which every order isomorphism is affine (in $f$). We refer to \cite{Lemmens2020,walsh:hal-02425988} for some technical conditions ensuring this. They unfortunately are not sufficient on spaces like that of continuous functions if not paired with homogeneity assumptions on $J$, in line with the work of \cite{Schffer1979}.

    \paragraph{\bf Structure of the article.} In \Cref{sec:preliminaries} we introduce the notions of irreducibility and of isomorphisms that we consider. \Cref{sec:general_sets} pertains to our main results on order isomorphisms, in particular when commuting with constants, with \Cref{sec:c-convex} devoted to $c$-convex functions. In \Cref{sec:examples}, we detail the consequences of these results for continuous and lower semicontinuous functions (\Cref{sec:examples_continuous}), for 1-Lipschitz functions (\Cref{sec:Busemann_Lip}), and finally for convex lower semicontinuous functions (\Cref{sec:examples_convex}).%
    
    \section{Preliminaries}\label{sec:preliminaries}
    \tb{Notation:}  The extended real line is denoted by $\pmR=[-\infty,+\infty]$, we also use the notations $\pR=(-\infty,+\infty]$, $\mR=[-\infty,+\infty)$ and $\R_{>0}=(0,+\infty)$. Given any sets $\X$ and $E$, we denote by $E^\X$ the set of functions $f$ from $\X$ to $E$ and $\Id_\X:\X\to\X$ the identity map. When $\X$ is equipped with a topology, a function $f:\X\rightarrow \pmR$ is said to be lower semicontinuous (l.s.c.) if its epigraph is a closed subset of $\X\times\pmR$. Its domain is defined by $\Dom(f)=\{x\in\X \, | \, f(x)<+\infty\}$. A function $f$ is said to be proper if its domain is non-empty and if it does not take the value $-\infty$. For $\G\subset\pmRX$ and some $x\in\X$, if there exists a function $g\in\G$ such that $g(x)\in\R$, we say that $\G$ is \emph{proper} at $x$, and more simply that $\G$ is \emph{proper} if it is \emph{proper} at all $x\in\X$. When $\X$ and $\Y$ are topological vector spaces, the set $\X^*$ denotes the dual space of $\X$ of continuous linear forms from $\X$ to $\R$ and $\Lcal(\X,\Y)$ refers to the space of continuous linear operators from $\X$ to $\Y$. We write the duality product over $\X^*\times\X$ as $\bk{\cdot,\cdot}:\X^*\times\X\to\R$.\\ %
    
    From now on, $\X$ and $\Y$ are two given nonempty sets.%
    \begin{Definition}\label{def:complete_subspace} We say that a subset $\G$ of $\pmRX$ is sup-stable (resp.\ finitely sup-stable) if it is stable under arbitrary sups (resp.\ pairwise sups), and that it is a complete subspace of $\pmRX$ if it is stable under arbitrary sups and addition of scalars. For $\G\subset \pmRX$, we define its sup and inf-closures as 
    \begin{align*}
        \Gsup&:=\{\sup_{\alpha\in\A} h_\alpha \, | \, \A \text{ an index set}, \, \{h_\alpha\}_{\alpha\in\A}\subset \G\}, \\
        \bGinf&:=\{\inf_{\alpha\in\A} h_\alpha \, | \, \A \text{ an index set}, \, \{h_\alpha\}_{\alpha\in\A}\subset \G\}.
    \end{align*}
     The set $\G$ is sup-stable if and only if $\G=\Gsup$. We say that $\G$ is inf-stable if $\G=\bGinf$.
    \end{Definition}
    The notion of complete subspace refers to idempotent (or tropical) functional analysis, in which spaces of functions are equipped with the supremum operation and with an additive action of real scalars, see~\cite{litvinovetal,McEneaney2006}. 
    \begin{Definition}\label{def:inf-relatif}
        Let $\G$ be a partially ordered set. Given an arbitrary index set $\Asc$, fix a family $(g_\alpha)_{\alpha\in \Asc}\in\G^{\Asc}$. If it exists, we define its infimum relatively to $\G$ as
		\[ \inf\nolimits^\G_\alpha g_\alpha:=\max\{h\in \G \, |\, \forall\alpha\in \Asc, \, h\le g_{\alpha}\} .\]%
	The latter set does admit a greatest element if $\G$ is sup-stable. In particular, if the family $(g_\alpha)_{\alpha\in\Asc}$ consists
		of a single element $g$, then,
    \begin{equation}\label{eq:def_infG}
          \infG g := \max\{h \in \G\,|\, h\leq g \}.
      \end{equation}
      Notably, if $g\in\G$, then $g=\infG g$.
    \end{Definition}
    Relative infima underline that infima or suprema can be done either w.r.t.\ the absolute (entry-wise) supremum of $\pmRX$ or its relative counterpart in $\G$. 
    \begin{Definition}\label{def:irreducible_points} 
        Let $\G$ be a subset of $\pmRX$. We say that a point $f\in\G$ is sup-irreducible if, for all $g,h\in\G$,
        \begin{equation}
            f=\sup(g,h) \, \implies \, f=g \text{ or } f=h.
        \end{equation}
        Similarly a point $f\in\G$ is inf-irreducible if, for all $g,h\in\G$,
        \begin{equation}
            f=\inf(g,h) \, \implies \, f=g \text{ or } f=h.
        \end{equation}%
        We say that a point $f\in\G$ is $\G$-relatively-inf-irreducible if, for all $g,h\in\G$,
        \begin{equation}
            f=\infG(g,h) \, \implies \, f=g \text{ or } f=h.
        \end{equation}
    \end{Definition}
     The term ``irreducible'' can be understood by analogy with irreducible polynomials.
    \begin{Claim}\label{claim:inf-relatif}
            Let $\G$ be a sup-stable subset of $\pmRX$. For all $\tilde g,\tilde h\in\pmRX$, if $f=\inf(\tilde g,\tilde h)$ with $f\in\G$, then $f=\inf(\infG\tilde g,\infG\tilde h)$. The $\G$-relatively-inf-irreducible functions are thus also inf-irreducible.%
        \end{Claim}
        \begin{proof} Take $f=\inf(\tilde g,\tilde h)$ with $f\in\G$. Since $f\le \tilde g$, we have that $f\le \infG(\tilde g)$ by definition of $\infG$ in \eqref{eq:def_infG}. The same holds for $\tilde h$, hence $f\le \tilde f:=\inf(\infG\tilde g,\infG\tilde h)$ and $f\le \infG \tilde f \le \tilde f$. Moreover $\tilde f \le\tilde g$ and $\tilde f \le\tilde h$, so $\tilde f \le f=\inf(\tilde g,\tilde h)$ and $f=\tilde f$. If $\tilde g,\tilde h\in\G$ and $f$ is $\G$-relatively-inf-irreducible, then w.l.o.g\ $f=\infG\tilde g=\tilde g$ so $f$ is inf-irreducible.
        \end{proof}
    We define the indicator functions $\delta^\bot_{x},\delta^\top_{x}\in \pmRX$ as follows
	\begin{equation}\label{eq:def_Dirac}
		\delta^\bot_{x}\left(y\right):=\left\{\begin{array}{ll}
			0 & \text { if } y=x, \\
			-\infty & \text { otherwise,} \end{array}\right. \quad 
		\delta^\top_{x}\left(y\right):=\left\{\begin{array}{ll}
			0 & \text { if } y=x, \\
			+\infty & \text { otherwise.}
		\end{array}\right.
	\end{equation}
    We will often informally refer to them as Dirac masses. Notably, if $\delta^\top_{x}\in\G$, then it is clearly $\G$-relatively-inf-irreducible, and if $\delta^\bot_{x}\in\G$, then it is sup-irreducible.

    \begin{Definition}\label{def:isomorphisms} 
        A map $J:\F \rightarrow \G$ where $\F$ and $\G$ are partially ordered sets is said to be
        \begin{enumerate}[labelindent=0cm,leftmargin=*,topsep=0.1cm,partopsep=0cm,parsep=0.1cm,itemsep=0.1cm,label=\roman*)]
			\item an order isomorphism if it is invertible and if this map and its inverse are both order preserving, i.e.\ for all $f,g\in\F$
            \begin{equation}
                f\ge g \, \Leftrightarrow \, Jf \ge Jg;
            \end{equation}
            \item a max-isomorphism if it is invertible and if it commutes with pairwise suprema, i.e.\ $J(\sup(f,g))=\sup(Jf,Jg)$, assuming that $\G$ and $\F$ are sup-stable;
			\item a (max,+)-isomorphism if it is a max-isomorphism and if we have $J(f+\la)=Jf +\la$ for $\la \in\R$, assuming that $\G$ and $\F$ are complete subspaces of $\pmRX$;
            \item order reversing if for all $f,g\in\F$
            \begin{equation}
                f\ge g \, \implies \, Jf \le Jg;
            \end{equation}
            \item an order anti-isomorphism if it is invertible and if this map and its inverse are both order reversing;
            \item an anti-involution if $J:\G \rightarrow \G$, $JJ=\Id_\G$ and $J$ is order reversing.
		\end{enumerate}
    \end{Definition}

    Order isomorphisms are a more general notion than max-isomorphisms, but the two coincide when $\F$ and $\G$ are sup-stable.%
    \begin{Lemma}\label{lem:order_implies_max} Let $\F\subset\pmRY$ and $\G\subset\pmRX$ be sup-stable sets. Then the order isomorphisms  $J:\F\rightarrow \G$ commute with an arbitrary suprema, and thus coincide with max-isomorphisms.  Moreover, if $J:\F\rightarrow \G$ is a max-isomorphism, then so is $J^{-1}:\G\rightarrow \F$. Anti-isomorphisms $\bF:\F\to\G$ send relative infima of $\F$ onto suprema of $\G$, i.e.\ $\bF(\inf\nolimits^\F_{\alpha\in \Asc} f_\alpha) = \sup_{\alpha\in \Asc}\bF( f_\alpha)$ for any $\{f_\alpha\}_{\alpha\in \Asc}\subset\F$.
    \end{Lemma}

    \begin{proof} Take $(f_\alpha)_{\alpha\in\Asc}\subset \F^\Asc$ for $\Asc$ some arbitrary index set and let $J:\F\to\G$ be an order isomorphism. Then $J(\sup_\alpha f_\alpha)\ge \sup_\alpha Jf_\alpha$ since $J$ preserves the order. The same holds for $J^{-1}$ so
       \begin{equation}\label{eq:sup_arbitraires_isophi}
       \sup_\alpha f_\alpha = J^{-1}J(\sup_\alpha f_\alpha)\ge J^{-1}(\sup_\alpha J f_\alpha) \ge \sup_\alpha J^{-1}J f_\alpha=\sup_\alpha f_\alpha,
        \end{equation}
        hence $J(\sup_\alpha f_\alpha)=\sup_\alpha J f_\alpha$. Let $J$ be a max-isomorphism, take $f,g\in\G$, $f_0,g_0\in\F$ such that $f=Jf_0$ and $g=Jg_0$. Then
    \begin{equation*}
        J^{-1}(\sup(f,g))=J^{-1}(\sup(Jf_0,Jg_0))=J^{-1}J(\sup(f_0,g_0))=\sup(f_0,g_0)=\sup(J^{-1}f,J^{-1}g).
    \end{equation*}
        
        Similarly, let $\bF^{-1}:\G\to\F$ be the inverse of $\bF$. As $\bF$ is an order anti-isomorphism, it satisfies that,
		 as $\inf\nolimits^\F_{\alpha\in \Asc} f_\alpha \le f_\beta$ for any $\beta\in \Asc$,
		\begin{displaymath}
			\bF(\inf\nolimits^\F_{\alpha\in \Asc} f_\alpha)\ge \sup_{\alpha\in \Asc}\bF( f_\alpha) \ge \bF( f_\beta).
		\end{displaymath}
		Since $\G$ is sup-stable, $\sup_{\alpha\in \Asc}\bF( f_\alpha)\in\G$, so composing by $\bF^{-1}$ and using that $\bF^{-1}$ is order reversing, we derive that
		\begin{equation*}
			\bF^{-1}\bF(\inf\nolimits^\F_{\alpha\in \Asc} f_\alpha)\le \bF^{-1}(\sup_{\alpha\in \Asc}\bF( f_\alpha)) \le \bF^{-1}\bF( f_\beta).
		\end{equation*}
		As $\bF$ is an anti-isomorphism, taking the infimum in $\F$ over $\beta$ on the r.h.s., yields
		\begin{equation*}
			\inf\nolimits^\F_{\alpha\in \Asc} f_\alpha  \le \bF^{-1}(\sup_{\alpha\in \Asc}\bF( f_\alpha)) \le \inf\nolimits^\F_{\beta\in \Asc} f_\beta
			.
		\end{equation*}
		So $\bF(\inf\nolimits^\F_{\alpha\in \Asc} f_\alpha) =  \sup_{\alpha\in \Asc}\bF( f_\alpha)$.%

    \end{proof}
  We will be sometimes interested in specific mappings in order to provide examples.
    \begin{Definition}\label{def:rmax-sesqui-lin}
		A map $B:\pmRY\rightarrow \pmRX$ is said to be $\pmRmax$-sesquilinear if $B(\inf\{f_i\}_{i\in I})=\sup\{Bf_i\}_{i\in I}$ and $B(f + \lambda)=Bf - \lambda$ (with $+\infty$ absorbing on the l.h.s.\, i.e.\ $f+\infty=+\infty$, and  $-\infty$ absorbing on the r.h.s., i.e.\ $Bf-\infty=-\infty$), for any finite index set $I$ and $\lambda\in\pmR$; we say in addition that $B$ is continuous if $B(\inf\{f_i\}_{i\in I})=\sup\{Bf_i\}_{i\in I}$ holds even for infinite families.
   
		The range $\Rg(B)$ of $B$ is defined as the set of functions $g\in \pmRX$ such that $g=Bf$ for some $f\in \pmRY$.
	\end{Definition}
    
	\begin{Proposition}[Theorem 3.1, \cite{singer1984conj}]\label{thm:singer}
		A map $\bB:\pmRY\rightarrow \pmRX$ is $\pmRmax$-sesquilinear and continuous if and only if there exists a kernel $b:\X\times \Y \rightarrow \pmR$ such that $\bB f(x)= \sup_{y\in\Y} b(x,y) - f(y)$ for all $f\in\pmR^{\Y}$ (with $-\infty$ absorbing). Moreover in this case $b$ is uniquely determined by $\bB$ as $b(\cdot,y)=\bB\delta^\top_{y}$.
	\end{Proposition}
     Hence the relation between $\pmRmax$-sesquilinear and continuous maps, a.k.a.\ \emph{(Fenchel-Moreau) conjugations}, and kernels is one-to-one, so it is equivalent to study a map or its kernel. Moreover \Cref{thm:singer} allows to write
	\begin{equation}\label{eq:RgB_tpsd}
		\Rg(B)=\{\sup_{y\in \Y} a_y+ b(\cdot, y) \,|\, a_y\in \mR \}.
	\end{equation}
    where $b:\X\times \Y \rightarrow \pmR$ is the kernel associated with $\bB$. Notice also that any complete subspace $\G$ of $\pmRX$ is a range for the specific choice $\Y=\G$ and $b(\cdot, g)=g(\cdot)$, since
     \begin{equation*}
         \G=\{\sup_{g\in \G} a_g+ g(\cdot) \,|\, a_g\in \mR \}.%
     \end{equation*}
     Obviously there may exist some more parcimonious representatives through other $b$ and $\Y$. Given a kernel $b:\X\times \Y \to \pmR$, let $\bB:\pmRY\rightarrow \pmRX$ and its transpose $\bB^\circ:\pmRX\rightarrow \pmRY$  be defined by 
    \begin{equation}
        \bB f(\cdot):= \sup_{y\in\Y} b(\cdot,y) - f(y), \quad \bB^\circ h(\cdot):= \sup_{x\in\X} b(x,\cdot) - h(x),\, \forall f\in\pmRY,\, h\in\pmRX
    \end{equation}
     with the convention that $-\infty$ is absorbing. Interesting examples of sets $\G=\Rg(B)$ are that of lower semicontinuous $C$-semiconvex functions for $b(x,y)=-\frac{C}{2}\|x-y\|^2$ for some  $C>0$, or of 1-Lipschitz functions for $b(x,y)=-d(x,y)$. We use a terminology rooted in abstract convexity theory, however these functions are nothing else than the $c$-convex functions in optimal transport theory, taking $c(x,y)=-b(x,y)$. The key relation is that $\bB=\bB \bB^\circ \bB$, see e.g.\  \cite[p.3, Theorem 2.1, Example 2.8]{Akian04setcoverings}. We provide now a characterization relating ranges and anti-isomorphisms.%
	\begin{Theorem}[Kernels and anti-isomorphisms]\label{thm:sym_anti-involution}
		Let $\G$ (resp.\ $\F$) be a complete subspace of $\pmRX$ (resp.\ $\pmRY$). Then the following statements are equivalent:
		\begin{enumerate}[label=\roman*),labelindent=0cm,leftmargin=*,topsep=0.1cm,partopsep=0cm,parsep=0.1cm,itemsep=0.1cm]
			\item \label{it_rgB} there exists a kernel $b:\X\times \Y \rightarrow \pmR$ such that $\G=\Rg(B)$ and $\F=\Rg(B^\circ)$;
			\item \label{it_antiinvol} there exists an order anti-isomorphism $\bF:\F\rightarrow \G$ commuting with the addition of scalars, i.e.\ $\bF(f + \lambda)=\bF f - \lambda$, for any $\lambda\in\R$ and $f\in\F$.
		\end{enumerate}%
		If these properties hold, then $\bF$ can be taken as the restriction of $\bB:\pmRY\to\pmRX$ to $\Rg(B^\circ)$. Moreover, for $\X=\Y$,  there exists an anti-involution $\bF$ over $\G$ if and only if there exists a symmetric $b$ such that $\G=\Rg(B)$.
	\end{Theorem}
     This result extends \cite[Theorem 5.1]{aubin2022tropical} originally written for anti-involutions only. The proof of \Cref{thm:sym_anti-involution} can be found in the Appendix on p\pageref{proof:sym_anti}. \Cref{thm:sym_anti-involution} exclusively provides the existence of the anti-isomorphism and does not discuss its uniqueness, which in general does not hold. It is one of the goals of this article to show that in many cases the map $\bF$ is nevertheless unique up to a translation $g$ and reparametrization $\phi$. 
    Anti-involutions were described by \cite[Definition 11]{ArtsteinAvidan2009} as a ``concept of duality'', this point of view being developed more at length in \cite{artstein2007}. \Cref{thm:sym_anti-involution} can also be seen as a functional analogue to \cite[Theorem 1.6]{ArtsteinAvidan2023} written for the power set $2^\X$. If there are no anti-involutions over a set $\G$, as that is the case for $\G$ the set of lower semicontinuous functions on a Hausdorff space with no isolated point, as shown in \Cref{lem:lsc_no_anti}, then there is no such duality for $\G$. Describing anti-isomorphisms over complete lattices as an abstract nonlinear duality is actually at the core of abstract convexity theory as summarized in \cite{singer1997}. In this article however, we do not use complete lattices since our sets are not inf-stable, unlike $\pmRX$, and we focus less on dualities and more on order isomorphisms.

    \section{Results on general sets $\G$}\label{sec:general_sets}

    \subsection{The shape of order (anti)isomorphisms}
    We start with a simple result showing that, if we know one reference anti-isomorphism, as given for instance by \Cref{thm:sym_anti-involution}, it is enough to study the order isomorphisms rather than the more arduous order anti-isomorphisms or anti-involutions.
    \begin{Lemma}\label{prop:isophi_anti-invol}
		Let $A,B:\F \rightarrow \G$ be two order anti-isomorphisms over partially ordered sets $\F$ and $\G$, and set $J=A^{-1}B$. Then $J$ is an order isomorphism over $\F$, and we have that $B=AJ$ and $A=BJ^{-1}$. In particular, if there exists an anti-involution $\bF:\G \rightarrow \G$, every anti-involution over $\G$ reads as $\bF J$ with $J:\G \rightarrow \G$ an order isomorphism satisfying $\bF J \bF J=\Id_{\G}$. The result also holds when ``order'' is replaced with (max,+).
	\end{Lemma}
	\begin{proof} 
        The operator $J$ has an inverse $J^{-1}=B^{-1}A$ and $J$ and $J^{-1}$ are order preserving, so $J$ is an order isomorphism. Moreover $\bF J$ is an anti-involution over $\G$ if and only if $\bF J \bF J=\Id_{\G}$.
    
	\end{proof}
    \begin{Remark}
        If $\G=\Rg(B)$ as in \eqref{eq:RgB_tpsd}, then, for every $f\in \Rg(B)$, there exists $(w_y)_{y\in\Y}$ such that $f(\cdot)=\sup_{y\in \Y} w_y+ b(\cdot, y)$. Hence we have that $J(f)=\sup_{y\in \Y} J(w_y+ b(\cdot, y))(\cdot)$ as $J$ commutes with the $\sup$, as shown in \eqref{eq:sup_arbitraires_isophi}. Thus to characterize $J$, we just have to study its action on the functions $\{w+b(\cdot,y)\}_{y\in\Y, w\in\mR}$. More generally this argument holds for sup-generators. However, unless there are other assumptions, such as continuity or compactness, see \Cref{lem:isophi_extreme_pts_generators_compact}, the $\{w+b(\cdot,y)\}_{y\in\Y, w\in\mR}$ are not necessarily sup-irreducible or relatively-inf-irreducible elements, which have the nice property of being stable under order isomorphisms.
    \end{Remark}%
    Let us now state how irreducible elements are invariants of order isomorphisms.
    \begin{Lemma}\label{lem:isophi_extreme_pts}
		Every order isomorphism $J:\F\rightarrow\G$ over sup-stable sets $\F\subset\pmRY$ and $\G\subset\pmRX$ sends sup-irreducible elements (resp.\ $\F$-relatively-inf-irreducible) of $\F$ onto sup-irreducible elements (resp.\ $\G$-relatively-inf-irreducible) of $\G$.

        Every anti-isomorphism $T:\F\to\G$ sends $\F$-relatively-inf-irreducible elements of $\F$ onto sup-irreducible elements of $\G$.\
	\end{Lemma}
	\begin{proof} 
		We just apply \Cref{lem:order_implies_max}. Let $f_0\in \F$ be sup-irreducible in $\F$. Take $g_1,h_1\in \G$ such that $J f_0=\sup(g_1,h_1)$. Since $J$ is onto, we can find $g_0,h_0\in \F$ such that $g_1=J g_0$, $h_1=J h_0$. Hence, using that $J$ commutes with the $\sup$, $J f_0=J(\sup(g_0,h_0))$. Since $J$ is bijective and $f_0$ is sup-irreducible, we have that $f_0=g_0$ or $f_0=h_0$, applying $J$, $J f_0$ is thus sup-irreducible in $\G$.
      For $\F$-relatively-inf-irreducible functions $f_0$ we have similarly that if $J f_0=\infG(g_1,h_1)$, then $ f_0=\sup \{ J^{-1}h \, | \, h\in\G,\, h\le g_1,\, h\le h_1\}$, so $f_0=\infnl^\F(J^{-1}g_1,J^{-1}h_1)$. 
      
      Concerning $T:\F\to\G$, if $Tf=\sup(g_1,h_1)$, then, by \Cref{lem:order_implies_max}, $f=\infnl^\F(T^{-1}g_1,T^{-1}h_1)$, and we conclude as before. %
	\end{proof}

    \begin{Remark}[Another notion of extremality based on comparability]\label{rmk:comparable_functions}
    In \cite[Definitions 4.14,4.15]{Artstein_Avidan_2012}, the authors say that a function $f\in\G$ is $P$-extremal if, for all $g,h\in\G$, $(f\le g) \& (f\le h)\Rightarrow g\le h \text{ or } h\ge g$, and conversely  $Q$-extremal if, for all $g,h\in\G$, $(f\ge g) \& (f\ge h)\Rightarrow g\le h \text{ or } h\ge g$. In other words the upper/lower bounds on an extremal $f$ should be ``comparable''. With the same proof as that of \Cref{lem:isophi_extreme_pts}, one can show that $P$ and $Q$-extremal functions are preserved by order isomorphisms. Moreover, if $f$ is $P$-extremal, then $f=\infG(g,h)\Rightarrow f\le g,h \stackrel{w.l.o.g}{\Rightarrow} g \le h \Rightarrow f=g$, so $f$ is relatively-inf-irreducible. For some sets, such as that of lower semicontinuous convex functions on a vector space, $P$-extremal functions and relatively-inf-irreducible ones coincide, since they are the ``Dirac masses'' $\delta^\top_x+\la$. However, for $1$-Lipschitz functions, there are simply no $P$-extremal functions,\footnote{This is shown for instance by taking two points $x_0\neq x_1$. Then any $1$-Lipschitz function $f$ satisfies $f\le \inf(d(\cdot,x_0)+\la_0, d(\cdot,x_1)+\la_0)$ with $\la_0=\max(f(x_0),f(x_1))$, but the two upper bounds are not comparable.} whereas the distance functions $d(\cdot,x_0)+\la$ are all relatively-inf-irreducible. This is why we favor the concept of irreducibility over comparability, even though the latter proved fruitful in \cite[Theorem 18]{ArtsteinAvidan2011hidden} or  \cite[Section 4]{Artstein_Avidan_2012}.
    \end{Remark}
 
    Sometimes, we may want to study sets that are not sup-stable, such as that of continuous functions or of Lipschitz functions, both of which are pointwise dense in the sup-stable set of lower semicontinuous functions. Under some regularity assumptions, we can however extend order isomorphisms over a set to the smallest sup-stable set containing it. The proof of the following is inspired by the theory of Galois connections, see \cite[Chapter 0, Section 3]{continuous}.
    \begin{Lemma}\label{lem:extension_order}
         Let $J:\F\rightarrow \G$ be an order isomorphism between some non-empty sets $\F\subset \pmR^\Y$ and $\G\subset \pmR^\X$ stable by addition of scalars and finitely sup-stable. Assume also that
          \begin{equation}\label{eq:h_halpha_majoration}
             \forall g\in\G, \, \forall g'=\lim_{\alpha} \uparrow g_\alpha\in\Gsup \text{ with $g_\alpha\in\G$, } (g\le g') \implies (\forall \epsilon>0,\, \exists \alpha_\epsilon \text{ s.t.\ } g\le g_{\alpha_\epsilon}+\epsilon),
         \end{equation}
         where $\lim_{\alpha} \uparrow g_\alpha$ denotes a pointwise non-decreasing net, and suppose that
         \begin{multline}\label{eq:J_continuity}
             \forall g'\in\Gsup, \,  \forall y \in \Y, \, \exists v_{y}:(0,+\infty)\to\R \text{ s.t.\ } \lim_{\epsilon'\to 0} v_y(\epsilon')=0 \text{ and }\\
             \forall g\in \G, (g\le g')\implies \left(\forall \epsilon>0,\, J^{-1}(g+\epsilon)(y)\le J^{-1}g(y)+v_y(\epsilon) \right).
         \end{multline}
         Then $J$ can be extended to an order isomorphism between $\Fsup$ and $\Gsup$.%

        Moreover, if for $f_0\in\R^\Y$ and $g_0\in\R^\X$, we have $\F'=f_0+\F$ and $\G'=g_0+\G$, then every order isomorphism $J':\F'\rightarrow\G'$ is of the form
        \begin{equation}\label{eq:isophi_translation}
            J'(f)=J(f-f_0)+g_0
        \end{equation}
        where $J:\F\rightarrow\G$ is an order isomorphism. The irreducible elements of $\F'$ are simply those of $\F$ translated by $f_0$. The results above also hold when ``order'' is replaced with (max,+), in particular \eqref{eq:J_continuity} is then automatically satisfied with $v_y(\epsilon)=\epsilon$.
    \end{Lemma}

    Eq.\eqref{eq:h_halpha_majoration} has to do with the regularity of the functions in $\G$, whereas \eqref{eq:J_continuity} is a form of uniform continuity of $J$.
    
    For instance, if $J^{-1}$ is $K$-Lipschitz for some $K>0$ in the sense $|J^{-1}(g+\epsilon)(y)-J^{-1}(g)(y)|\le K \epsilon$ for all $\epsilon>0$, $g\in\G$, $y\in\Y$, then \eqref{eq:J_continuity} holds. On the other hand, \eqref{eq:h_halpha_majoration} is satisfied if we have that $\G$ is a subset of the continuous functions on a compact set $\X$, as a consequence of Dini's theorem applied to $g'_\alpha=\sup(g,g_\alpha)$. Eq.\eqref{eq:h_halpha_majoration} still holds if $\X$ is no longer compact but the functions ``tend to $0$ at infinity'', i.e.\ for all $h\in\G$, $\epsilon>0$, there exists a compact $K_\epsilon\subset \X$ s.t.\ for $x\notin K_\epsilon$, $|h(x)|\le\epsilon$.

	\begin{proof} %
    Define the two candidate extensions $\bJ:\Fsup\rightarrow\Gsup$ and $\bJ^1: \Gsup\rightarrow\Fsup$ of $J$ and $J^{-1}$%
	    \begin{align}
	        \bJ f &= \sup \{ J(h) \, | \, h\in \F, \, h\le f \},\\
         \bJ^1 g &= \sup \{ J^{-1}(g') \, | \, g'\in \G, \, g'\le g \}.
	    \end{align}
    Clearly their restrictions to $\F$ and $\G$ are equal to $J$ and $J^{-1}$ respectively. Moreover $\bJ$ and $\bJ^1$ are order preserving. Fix $f\in \Fsup$, and set $\tilde f=\bJ^1 \bJ f$. We want to show that $f=\tilde f$. We have that
    \begin{equation*}
        \tilde f=\bJ^1 \bJ f= \sup \{ J^{-1}(g) \, | \, g\in \G, \, g\le \sup_{h'\in\F, h'\le f} J(h') \}=\sup \{ h \, | \, h\in \F, \, J(h)\le \sup_{h'\in\F, h'\le f} J(h') \}.
    \end{equation*}
     Let us show that $\tilde f \le f$. Fix $h\in \F$  such that $ J(h)\le \sup_{h'\le f} J(h')$. Set $g=J(h)$ and $g'=\sup_{h'\le f} J(h')$. Since $\G$ is finitely sup-stable, we can take $g'=\lim\uparrow g_\alpha$ with each $g_\alpha$ a supremum of a finite number of functions $J(h')$ satisfying $h'\le f$. Hence by \eqref{eq:h_halpha_majoration}, for all $\epsilon>0$, there exists $g_{\alpha_\epsilon}$ such that $g\le g_{\alpha_\epsilon}+\epsilon$, and $g_{\alpha_\epsilon}=\sup_{h_i'\le f, \, i\in\{1,\dots, N\}} J(h_i')$. Notice that in \Cref{lem:order_implies_max}, eq.\eqref{eq:sup_arbitraires_isophi} still holds in our setting for finite $\Asc$. Hence $J^{-1}$ commutes with a finite number of suprema. So we obtain
     \begin{equation*}
         h=J^{-1}(g)\le J^{-1}(g_{\alpha_\epsilon}+\epsilon)= \sup_i J^{-1}(J(h'_i) +\epsilon).
     \end{equation*}
     Applying \eqref{eq:J_continuity}, since $J(h'_i)\le g'$, for all $y\in\Y$, it holds $h(y)\le \sup_i h'_i(y)+v_{y}(\epsilon)\le  f(y)+ v_{y}(\epsilon)$. Taking $\epsilon$ to $0$, we get $h \le f$. Taking the supremum over $h$ shows that $\tilde f \le f$.
     
    To prove that $\tilde f \ge f$, we are going to show that, for any $f\in\Fsup$ and $g\in \Gsup$, we have that $g\ge \bJ f \implies \bJ^1 g \ge f$. Take any $f,g$ satisfying $g\ge \bJ f$. Fix $f'\in \F$ such that $f'\le f$, so $Jf' \le \bJ f \le g$, whence $f'\le \bJ^1 g$. Taking the supremum over $f'$, we obtain $\bJ^1 g \ge f$ as claimed. Since $\bJ f \ge \bJ f$, the implication we have just shown gives that $\tilde f=\bJ^1 \bJ f \ge f$. We thus obtained that $f=\bJ^1 \bJ f$, so $\bJ$ is an order isomorphism.%

    For \eqref{eq:isophi_translation}, just set $J(f):=J'(f+f_0)-g_0$ with inverse $J^{-1}(g)=(J')^{-1}(g+g_0)-f_0$. Take a sup-irreducible $f\in\F$, then, for all $g,h\in\F$, $f+f_0=\sup(g+f_0,h+f_0)\implies f=g \text{ or } f=h$, so $(f+f_0)$ is sup-irreducible in $\F'$. The other implications for irreducible functions are shown similarly.
	\end{proof}
 
    It was rightly mentioned in \cite[p670]{ArtsteinAvidan2009} that most results could be obtained either starting from order preserving or order reversing operators, without much difference. However we will focus on order preserving ones, as done in the convex case by \cite[p75]{Iusem2015}. Indeed order preserving operators have the advantage that the identity operator is always order preserving, whereas we will meet spaces with no anti-involution. Most of the following results stem from characterizing either the sup or inf-irreducible elements. The usefulness of focusing on such elements was already exemplified in \cite[p5]{GK} for finite sets $\X$.%

    \subsection{The $e_x$, the Dirac-like inf-irreducible functions of $\G$}
    
    We now turn to our candidates for inf-irreducible elements of a complete subspace $\G$ of $\pmRX$. Define the Archimedean class of a function $f\in\G$ as
    \begin{equation}
        [f]:=\{g\in\G \, | \, \exists \alpha \, \in\R,\, f-\alpha \le g \le f+\alpha\}.
    \end{equation}
    Let us put an order on Archimedean classes, saying that $[f]\le [g]$ if there exists $\alpha \in \R$ such that $f\le g+\alpha$. We say that $[f]$ is maximal if  $[f]\le [g] \implies [f]\ge [g]$.
    \begin{Claim}\label{claim:archimd_preserved}
        Let $\G$ (resp.\ $\F$) be a complete subspace  of $\pmRX$ (resp.\ $\pmRY$) for some sets $\X$ and $\Y$. Let $J$ be a (max,+)-isomorphism from $\F$ onto $\G$. If $f\in\F$ is such that $[f]$ is maximal, then $[Jf]$ is also maximal.
    \end{Claim}
    \begin{proof}
        Take $g_0\in\G$ and $\alpha\in\R$ such that $Jf\le g_0+\alpha$. Fix $f_0\in\F$ such that $g_0=Jf_0$, then $f\le f_0+\alpha$. Since $[f]$ is maximal, there exists $\alpha'\in\R$, $f\ge f_0+\alpha'$, composing with $J$, we obtain that $Jf\ge g_0+\alpha'$, so $[Jf]=[g_0]$ and $[Jf]$ is maximal.
    \end{proof}
    It is thus tempting to focus on maximal Archimedean classes, and we will do so in \Cref{thm:max_plus_isophi_general_ss_hyp_simple}. However some sets have no maximal classes as seen in Remark~\ref{rmk:no_max_archimedean}.
    
    We will often need in the sequel to consider at least some finite values, hence we restrict ourselves to proper sets $\G$ of functions of $\pmRX$. We are also not interested in situations where points cannot be distinguished by $\G$, since we could then quotient $\X$ by a $\G$-induced equivalence relation. Hence we introduce
    \begin{Assumption}\label{assumption:proper} The set $\G\subset\pmRX$ is both proper and point separating, i.e.\ for any $x, x'\in\X$ with $x\neq x'$, there exists $g_1,g_2\in\G$ such that $g_1(x),g_2(x),g_1(x'),g_2(x')\in\R$ and $g_1(x)-g_1(x')\neq g_2(x)-g_2(x')$.
    \end{Assumption}
    
    \begin{Lemma}\label{lem:prop_e_x} 
    Let $\G$ be a complete subspace of $\pmRX$ for some set $\X$, with $\G$ satisfying \Cref{assumption:proper}. Define, for any $x\in\X$, the function $e_x:\X\rightarrow \pmR$ by
    \begin{equation}\label{eq:def_e_x}
        e_x(\cdot):=\sup \{u \in \G \,|\, u(x)\le 0\}.
    \end{equation}%
     Then $e_x\in\G$, $e_x(x)=0$, $e_x$ is inf-irreducible in $\G$ and for all $x\neq x'$ and $\la\in\pmR$, $e_x\neq e_{x'}+\la$. We also have 
     \begin{equation}\label{eq:e_x_formule}
         e_x(y)=\sup_{u\in\G}u(y)-u(x),
     \end{equation} 
     with $-\infty$ absorbing. Moreover, for any $f\in\G$, we have the representation (with $+\infty$ absorbing)
    \begin{equation}\label{eq:representation_formula_e_x}
        f=\inf_{x\in\Dom(f)} e_x+f(x), 
    \end{equation}
     and $e_x$ satisfies a triangle inequality: for all $x,x',x''\in\X$, $e_x(x'')\le e_{x'}(x'')+e_{x}(x')$. If $f\in\G$ is such that $[f]$ is maximal, then for all $x_0\in\Dom(f)$ we have $ [e_{x_0}]= [f]$, i.e.\ we can fix $\la_0\in\R$, such that $e_{x_0}+\la_0 \le f$. If $\G$ contains the constant functions, then $e_x$ has nonnegative values.
    \end{Lemma}
    \begin{proof}%
        Let $x\in\X$. Since $\G$ is proper, we can fix $u_x\in\G$ such that $u_x(x)\in\R$. So $e_x\ge u_x-u_x(x)$, whence $e_x(x)\ge 0$, while by definition $e_x(x)\le 0$. We have that $e_x\in\G$ as $\G$ is sup-stable. If $e_x=\inf(f,g)$ for some $f,g\in\G$ then w.l.o.g.\ we can assume that $f(x)=e_x(x)=0$ so, by maximality of $e_x$, $f\le e_x=\inf(f,g)$ and $f=e_x$. Take $f\in\G$ with $[f]$ maximal. Since $f\le e_x+f(x)$ for all $x\in\Dom(f)$ i.e.\ $[e_x]\ge [f] $, we have that $[e_x]\le [f] $ implies that $[e_x]= [f]$. Fix $x\neq x'$ and $\la\in\pmR$ and $g_1,g_2\in\G$ as in \Cref{assumption:proper}. Assume that $e_x= e_{x'}+\la$. We can already exclude $\la\in\{\pm \infty\}$ since $e_x$ has a finite value in $x$. Moreover by definition of $e_x$ we have
        \begin{equation*}
            e_{x'}+\la=e_{x}\ge g_1-g_1(x), \quad e_{x}-\la=e_{x'}\ge g_1-g_1(x').
        \end{equation*}
        Evaluating in $x'$ and $x$ the above expressions gives $\la\ge g_1(x')-g_1(x)\ge \la$. The same holds for $g_2$, which contradicts \Cref{assumption:proper}. So $e_x\neq e_{x'}+\la$. Eq.\eqref{eq:e_x_formule} is just a rewriting of \eqref{eq:def_e_x}.

        Concerning \eqref{eq:representation_formula_e_x}, on the one hand, we can take the infimum over $x'\in\X$ of $f\le e_{x'}+f(x')$, the latter giving the triangle inequality when taking $f=e_{x}$. On the other hand, for any $y\in\X$, we have that $\inf_{x\in\Dom(f)} e_x(y)+f(x)\le e_y(y)+f(y)=f(y)$.
    \end{proof}
    
    \begin{Remark}[Increasing functions have no maximal Archimedean class]\label{rmk:no_max_archimedean} Consider the set $\G$ of increasing functions from $\R$ to $\pmR$. It is sup and inf-stable. A quick computation gives $e_x=\delta^\top_{\{ y\in\R | y\le x\}}$ where $\delta^\top_A$ denotes the indicator function of a subset $A\subset \R$. This $e_x$ looks like a Dirac function, but $[e_x]$ is strictly smaller than $[e_{x'}]$ for any $x'< x$, so no Archimedean class is maximal. Moreover the constant functions are inf-irreducible, but do not belong to $\{e_x+\la\}_{\la\in\R,x\in\R}$.
    
    \end{Remark}
    
    Even though $e_x$ is defined through $\G$, it is actually deeply related to $\bGinf$, for which it is an ``inf-generator''. %
    \begin{Lemma}[The $e_x$ and the inf-closure $\bGinf$]\label{lem:e_x_and_Ginf}
       Consider the inf-closure $\bGinf$ of a complete subspace $\G\subset\pmRX$  satisfying \Cref{assumption:proper}. Then any $f\in\G$ such that $f$ is inf-irreducible in $\bGinf$ is inf-irreducible in $\G$ and
        \begin{equation}\label{eq:e_x_Ginf}
            \bGinf =\{  \inf_x e_x(\cdot)+w_x \, | \, w_x \in \pmR \} \text{   and    } e_x(\cdot)=\sup \{ u \in \bGinf \, | \, u(x)<0\}.
        \end{equation}        
    \end{Lemma}
    \begin{proof}
        The first statement is just a consequence of the fact that $\G\subset \bGinf$. Concerning \eqref{eq:e_x_Ginf}, since $e_x(\cdot)+w_x\in\G$, by definition $\{  \inf_x e_x(\cdot)+w_x \, | \, w_x \in \pmR \}\subset \bGinf$. Take $f\in \bGinf$, so $f=\inf_{\alpha\in\A} h_\alpha$ for $ \A$ an index set and $(h_\alpha)_{\alpha\in\A}\subset \G$. By \eqref{eq:representation_formula_e_x}, each $h_\alpha$ has a representation by the $e_x$, commute the infima and set $w_x=\inf_{\alpha\in\A} h_\alpha(x)$ to conclude that $ \bGinf =\{  \inf_x e_x(\cdot)+w_x \, | \, w_x \in \pmR \}$.

        Set $\tilde e_x^<(\cdot)=\sup \{ u \in \bGinf \, | \, u(x)<0\}$ and $\tilde e_x(\cdot)=\sup \{ u \in \bGinf \, | \, u(x)\le 0\}$. By definition, $ \tilde e_x^<\le \tilde e_x$, and, for all $\epsilon>0$, $ \tilde e_x-\epsilon\le \tilde e_x^<$. Taking the limit $\epsilon\to 0^+$, $ \tilde e_x^<= \tilde e_x$. Since $\G\subset \bGinf$, $ e_x \le \tilde e_x$. Let $u=\inf_{\alpha\in\A} h_\alpha \in \bGinf$ such that $u(x)< 0$ and fix $\alpha\in\A$ such that $h_\alpha(x)\le 0$. As $h_\alpha\in\G$, we have $u\le h_\alpha \le e_x$. Taking the supremum over such $u$ implies that $\tilde e_x^< \le e_x $, which concludes the proof.
    \end{proof}
    For sets $\G=\Rg(B)$ we simply have that $e_x$ is the image of
    $b(x,\cdot)$.
    \begin{Lemma}[Extreme points $e_x(\cdot)$ and kernels $b(\cdot,y)$]\label{rmk:e_x_et_b_x}
            Let $\G=\Rg(B)$ as in \eqref{eq:RgB_tpsd} then, for all $x\in\X$ and $y\in\Y$,
            \begin{align}\label{eq:e_x_and_b_x}
                \bB^\circ e_x (y)=b(x,y), \quad e_x=\bB b(x,\cdot)=\sup_{y\in\Y}b(\cdot,y)-b(x,y).
            \end{align}
    \end{Lemma}
    \begin{proof} Fix $x\in\X$ and $y\in\Y$, then
         \begin{align*}
                \bB^\circ e_x (y)=\sup_{z\in\X} [b(z,y) - \sup_{u\in\G} [u(z)-u(x)]]=\sup_{z\in\X} \inf_{u\in\G} b(z,y) -u(z)+u(x).
            \end{align*}
                Take first $u=b(\cdot,y)$ to show that $\bB^\circ  e_x (z)$ is smaller than $b(x,y)$ and then $z=x$ to show it is also larger. Hence $\bB^\circ  e_x (y)=b(x,y)$. And conversely since $\bB\bB^\circ\bB=\bB$, $\bB^\circ$ is a (max,+) anti-isomorphism over $\G$ and $e_x=\bB b(x,\cdot)$.
    \end{proof}
    
    Using \eqref{eq:e_x_formule}, taking $\G$ to be the set of 1-Lipschitz functions over a metric space $(\X,d)$ leads to $e_x=d(\cdot,x)$, and, more generally, for $(C,p)$-Hölder functions, $e_x=Cd(\cdot,x)^p$ for $C>0$ and $p>0$. When $\G$ is the set of (convex) lower semicontinuous functions over a topological space $\X$, we have $e_x=\delta_x^\top$. Actually, whenever $\{\delta_x^\top\}_{x\in\X}\subset\G$, the Dirac masses are necessarily the $e_x$. There are many spaces for which it is the case, for instance all those that are defined through a growth condition, see \cite[p127]{singer1997}. We now turn to our first main result.
    
     \begin{Theorem}[Main theorem on (max,+)-isomorphisms]\label{thm:max_plus_isophi_general}%
        Let $\G$ (resp.\ $\F$) be a complete subspace  of $\pmRX$ (resp.\ $\pmRY$) for some sets $\X$ and $\Y$, with both $\F$ and $\G$ satisfying \Cref{assumption:proper}. Let $J$ be a (max,+)-isomorphism from $\F$ onto $\G$ and define $e_x\in\G$ (resp.\ $e'_y\in\F$) as in \eqref{eq:def_e_x}. Then the following statements are equivalent:
		\begin{enumerate}[label=\roman*),labelindent=0cm,leftmargin=*,topsep=0.1cm,partopsep=0cm,parsep=0.1cm,itemsep=0.1cm]
			\item \label{it_J_e_x} for all $y \in \Y$, $\la'\in\pmR$,  $J(e'_y+\la')\in\{e_x+\la\}_{x\in\X,\,\la\in\pmR}$;%
			\item \label{it_g_f_phi} there exists a function $g:\X\rightarrow \R$ and a bijection $\phi:\X \rightarrow \Y$ such that
            \begin{equation}\label{eq:charac_J}
			Jf(x)=g(x)+f(\phi(x))
		\end{equation}
            and, for all $f\in\F$, $h\in\G$, $(g+f\circ\phi)\in\G$ and $(-g\circ\phi^{-1}+h\circ\phi^{-1})\in\F$;
        \item \label{it_inf_ex} for all $x\in \X$, we have
                 \begin{equation}\label{eq:charac_J_ex}
			0=\inf_{y\in\Y} J e'_y(x)+J^{-1}e_x(y)
		\end{equation}
            and the infimum is attained at a unique point $\phi(x)\in\Y$.
		\end{enumerate}
     Conversely every $J$ satisfying \ref{it_g_f_phi} is order preserving, commutes with the addition of scalars, maps $\F$ to $\G$, and has an inverse with the same properties, so it is a (max,+)-isomorphism. Moreover the representation $(g,\phi)$ in \eqref{eq:charac_J} is unique. 
     \end{Theorem}
     The representation~\eqref{eq:charac_J} is analogous to the polar decomposition of matrices in the general real linear group as a product of a symmetric positive matrix and an orthogonal matrix.
	\begin{proof} 
        \ref{it_J_e_x}$\Rightarrow$\ref{it_g_f_phi} Since $J$ is order preserving and bijective, the constant function $c_\infty:x\mapsto +\infty$ is mapped onto itself. Indeed $Jc_\infty \ge \sup_{f\in \F} Jf = \sup_{h\in \G} h = c_\infty$. The same holds for $-\infty$.
        
        Fix $y\in\Y$ and recall that $e'_y(\cdot):=\sup \{v \in \F \,|\, v(y)\le 0\}$. By assumption \ref{it_J_e_x}, $Je'_y=e_{\psi(y)}+\bar g(y)$ with $\bar g(y)=Je'_y(\psi(y))\in\R$. Conversely there exists a $\varphi:\X\to\Y$ and $\la_x\in\R$ such that $e_{x}=J(e'_{\varphi(x)}+\la_x)=e_{\psi(\varphi(x))}+\la_x+\bar g(\varphi(x))$. However from \Cref{lem:prop_e_x}, we know that \Cref{assumption:proper} then ensures that $x=\psi(\varphi(x))$, so $\psi$ is surjective. Similarly, consider $Je'_y=Je'_{y'}$, to get that $\psi$ is injective.

        Take any $f\in \F$ and $h\in\G$. Then $f \le e'_y +f(y)$ and $h \le e_x +h(x)$. So composing with $J$ and $J^{-1}$, the following holds
        \begin{align*}
            (Jf)(\psi(y)) &\le e_{\psi(y)}(\psi(y))+\bar g(y) + f(y)=\bar g(y) + f(y),\\
            (J^{-1}h)(\psi^{-1}(x)) &\le e'_{\psi^{-1}(x)}(\psi^{-1}(x)) -\bar g(\psi^{-1}(x)) + h(x)=-\bar g(\psi^{-1}(x)) + h(x).
        \end{align*}
        Let us choose $h=Jf$ and $y=\psi^{-1}(x)$, thus
        \begin{equation*}
            f(y)\le -\bar g(y) + (Jf)(\psi(y))\le -\bar g(y)+\bar g(y)+f(y).
        \end{equation*}
        Hence $(Jf)(\psi(y))= \bar g(y) + f(y)$. We just take $\phi=\psi^{-1}$ and $g=\bar g\circ \psi^{-1}$ to derive \eqref{eq:charac_J}.

        \ref{it_J_e_x}$\Leftarrow$\ref{it_g_f_phi} Let $J$ as in \eqref{eq:charac_J}, which clearly maps the functions $\pm\infty$ onto themselves and is a (max,+)-isomorphism. Fix $y\in\Y$ and $\la'\in\R$. Since $J$ commutes with $(\sup,+)$-operations,
    \begin{align*}
        J (e'_y+\la')&=\la'+\sup \{Jv \,|\,v \in \F,\, v(y)\le 0\}=\la'+\sup \{u \in \G \,|\, J^{-1}u(y)\le 0\}\\
        &=\la'+\sup \{u \in \G \,|\, -g(\phi^{-1}(y))+u(\phi^{-1}(y))\le 0\}\\
        &=\la'+g(\phi^{-1}(y))+\sup \{u \in \G \,|\, u(\phi^{-1}(y))\le 0\}=\la'+g(\phi^{-1}(y))+e_{\phi^{-1}(y)},
    \end{align*}
    so \ref{it_J_e_x} holds. Assume that $Jf=g+f\circ \phi=\tilde g+f\circ \tilde \phi$ for some alternative $\tilde g$ and $\tilde \phi$ satisfying \eqref{eq:charac_J}. Then the above computation gives that $g(\phi^{-1}(y))+e_{\phi^{-1}(y)}=\tilde g(\tilde \phi^{-1}(y))+e_{\tilde \phi^{-1}(y)}$ for all $y\in\Y$. \Cref{lem:prop_e_x} then gives that necessarily $\phi=\tilde \phi$, and, evaluating in $\phi^{-1}(y)$, $g=\tilde g$. So the representation \eqref{eq:charac_J} is unique.

        \ref{it_g_f_phi}$\Rightarrow$\ref{it_inf_ex} Fix $x\in\X$. Notice that by \eqref{eq:e_x_formule}
         \begin{align*}%
        \inf_{y} J e'_y(x)+J^{-1}e_x(y)&=\inf_{y}\sup_{u\in\F} Ju(x) -u(y)+\sup_{v\in\G} J^{-1}v(y) -v(x)\\
        &=\inf_{y}\sup_{u,v\in\F} Ju(x) -u(y)+v(y) -Jv(x)\stackrel{u=v}{\ge}0
        \end{align*}
        Taking $y=\phi(x)$ and using \eqref{eq:charac_J}, we obtain \eqref{eq:charac_J_ex}.

        \ref{it_g_f_phi}$\Leftarrow$\ref{it_inf_ex} Since \eqref{eq:charac_J_ex} reads
        \begin{equation*}
			0=\inf_{y\in\Y} J e'_y(x)+J^{-1}e_x(y)= \inf_{y\in\Y}\sup_{u,v\in\F} Ju(x) -u(y)+v(y) -Jv(x),
        \end{equation*}
        we see that, for $y=\phi(x)$, we have $0=\sup_{u,v\in\F} Ju(x) -u(y)+v(y) -Jv(x)$. Hence we obtain $Ju(x) -u(y)+v(y) -Jv(x)\le 0$ for all $u,v\in \F$. But we also have $Jv(x) -v(y)+u(y) -Ju(x)\le 0$. The sum of these nonpositive terms is null, so they are both null, i.e.\
    \begin{equation*}
        Ju(x) -u(\phi(x))+v(\phi(x)) -Jv(x)=0, \quad \forall u,v\in\F.
    \end{equation*}
    Taking $v=e'_{\phi(x)}$, we have thus $Ju(x)=Je'_{\phi(x)}(x)+u(\phi(x))$, hence \eqref{eq:charac_J} holds.
    	\end{proof}

    The equivalence between the statements indicates that an explicit structure such as \Cref{thm:max_plus_isophi_general}~\ref{it_g_f_phi} can be retrieved from ``sending $\F$-Dirac masses onto $\G$-Dirac masses'' as required by \Cref{thm:max_plus_isophi_general}~\ref{it_J_e_x} and that this is rooted in a symmetrical property between $J$ and $J^{-1}$ with some infimum being attained as in \Cref{thm:max_plus_isophi_general}~\ref{it_inf_ex}. However all three properties depend on $J$ and we would like a general characterization of (max,+)-isomorphisms, for instance showing that they are all of the form \eqref{eq:charac_J}. This clearly will require some assumptions on $\X,\Y,\G,\F$ as seen from Remark~\ref{rmk:no_extension} exhibiting an isomorphism not satisfying \eqref{eq:charac_J}.

    \begin{Remark}[A counter-example concerning  isomorphisms]\label{rmk:no_extension}
      In general not every (max,+)-isomorphism is of the form \eqref{eq:charac_J}, and even sometimes none such $J$ exists. Consider a finite case where $\X=\Y=\{1,2,3\}$, and $\X'=\Y'=\{1,2\}$, and consider the kernels $a$, $b$, $c$, $d$ given by the following matrices %
      \[
      A= \left(\begin{array}{ccc} 0 & -2 & 0 \\ -2 & 0 & 0 \\ 0 & 0 & 0 \end{array}\right),
      \quad
      B= \left(\begin{array}{ccc} 0 & -1 & 0 \\ -3 & 0 & 0 \\ 0 & 0 & 0 \end{array}\right),
      \quad
      C= \left(\begin{array}{cc} 0 & -2  \\ -2 & 0 \end{array}\right),
      \quad
       D= \left(\begin{array}{cc} 0 & -1  \\ -3 & 0 \end{array}\right).
       \]
       I.e., we identify the matrix $A=(a_{ij})$ to the kernel $a: \{1,2,3\}^2 \to \R\cup\{-\infty\}$,
       $(i,j)\mapsto a_{ij}$, and make similar identifications for the matrices $B,C,D$.
       We shall denote by $\Rg_\R(A)$ the ``finite'' part of $\Rg(A)$, i.e.,
       $\Rg_\R(A) = \Rg(A)\cap \R^3=\{\sup_{1\leq j\leq 3} A(\cdot,j)+ \lambda_j \mid \lambda _j \in \R\}$.
       Observe that $\Rg(A)$ is the union of $\Rg_\R(A)$ and of the two vectors that are either identically
       $-\infty$ or $+\infty$. Hence to check the isomorphism
       of spaces of the form $\Rg(A)$, we shall limit our attention to their ``finite parts''.
        
	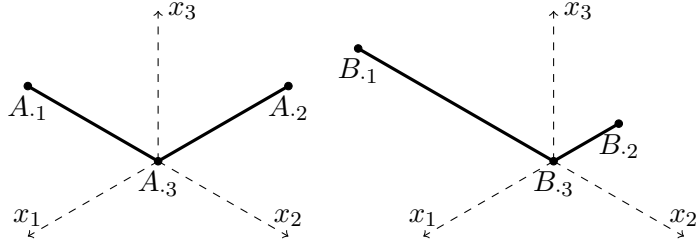
\begin{figure}[!ht]
		\begin{center}
			
			\def\coord#1#2#3{{-sqrt(3)/2*(#1-#2)} ,{ -(1/2)*#1 - (1/2)*#2 + #3}}
			
			\begin{tikzpicture}
				
				\coordinate (v1) at (\coord{0}{-2}{0});
				\coordinate (v2) at (\coord{-2}{0}{0});
				\coordinate (v3) at (\coord{-1}{0}{1});
				\coordinate (v4) at (\coord{0}{-1}{0});
				\coordinate (v5) at (\coord{0}{-0.5}{0});
				\coordinate (v6) at (\coord{1}{0}{0});
				\coordinate (v7) at (\coord{0}{0}{1});
                                				\coordinate (v0) at (\coord{0}{0}{0});

				\coordinate (p1) at (\coord{1}{0}{0});
				\coordinate (p2) at (\coord{0}{0}{1});
				\coordinate (p3) at (\coord{0}{-1}{0});

          				\draw[gray,draw=black,very thick] (v1) -- (v0) -- (v2);

				\draw[dashed,->] (\coord{0}{0}{0}) -- (\coord{2}{0}{0}) node[above] {$x_1$};
				\draw[dashed,->] (\coord{0}{0}{0}) -- (\coord{0}{2}{0}) node[above] {$x_2$};
				\draw[dashed,->] (\coord{0}{0}{0}) -- (\coord{0}{0}{2}) node[above, right] {$x_3$};
				
				\filldraw (v1) circle (0.3ex) node[below] {$A_{\cdot 1}$};
				\filldraw (v0) circle (0.3ex) node[below] {$A_{\cdot 3}$};
				\filldraw (v2) circle (0.3ex) node[below] {$A_{\cdot 2}$};
				
			\end{tikzpicture}
			\begin{tikzpicture}
				
				\coordinate (v1) at (\coord{0}{-3}{0});
				\coordinate (v2) at (\coord{-1}{0}{0});
				\coordinate (v3) at (\coord{-1}{0}{1});
				\coordinate (v4) at (\coord{0}{-1}{0});
				\coordinate (v5) at (\coord{0}{-0.5}{0});
				\coordinate (v6) at (\coord{1}{0}{0});
				\coordinate (v7) at (\coord{0}{0}{1});
                                \coordinate (v0) at (\coord{0}{0}{0});

				\coordinate (p1) at (\coord{1}{0}{0});
				\coordinate (p2) at (\coord{0}{0}{1});
				\coordinate (p3) at (\coord{0}{-1}{0});

          				\draw[gray,draw=black,very thick] (v1) -- (v0) -- (v2);

				\draw[dashed,->] (\coord{0}{0}{0}) -- (\coord{2}{0}{0}) node[above] {$x_1$};
				\draw[dashed,->] (\coord{0}{0}{0}) -- (\coord{0}{2}{0}) node[above] {$x_2$};
				\draw[dashed,->] (\coord{0}{0}{0}) -- (\coord{0}{0}{2}) node[above, right] {$x_3$};
				
				\filldraw (v1) circle (0.3ex) node[below] {$B_{\cdot 1}$};
				\filldraw (v0) circle (0.3ex) node[below] {$B_{\cdot 3}$};
				\filldraw (v2) circle (0.3ex) node[below] {$B_{\cdot 2}$};
				
			\end{tikzpicture}
		\end{center}
		\caption{The isomorphic spaces $\Rg_\R(A)$ and $\Rg_\R(B)$ (here intersected with the hyperplane $x_1+x_2+x_3=0$) do not have (max,+)-isomorphisms of the form \eqref{eq:charac_J}.}\label{fig-flag}

        \end{figure}
       First, it is immediate that $\Rg(C) $ and $\Rg(D)$ are isomorphic. Indeed, $\Rg_\R(C) =\{x\in \R^2\mid -2+x_1\leq x_2 \leq 2+x_1\}$
       and $\Rg_\R(D) = \{x\in \R^2\mid -3+x_1\leq x_2 \leq 1+x_1\}$, and the translation $(x_1,x_2)\mapsto(x_1,x_2-1)$
       is a (max,+) isomorphism sending $\Rg_\R(C)$ to $\Rg_\R(D)$.

       We also remark that an element in the interior of $\Rg_R C$ can be written in a unique
       way as a (max,+) linear combination of the columns of $C$. Using this observation, we deduce
       that the map which sends the $j$-th column of $C$ to the $j$-th column of $A$
       extends to an injective (max,+) linear map from $\Rg_\R(C)$ to $\Rg_\R(A)$.
       Moreover, the map from $\R^3$ to $\R^2$ which acts by restriction
       to the first coordinate is an inverse of this map. Observe also that the third
       column of $A$ is a (max,+) linear combination (actually, the supremum) of its first two columns.
       In that way, we deduce that $\Rg_\R(A)$ and $\Rg_\R(C)$ are isomorphic.
       A similar argument shows that $\Rg_\R(B)$ and $\Rg_\R(D)$ are isomorphic.
       Using $\Rg_\R(C)\simeq \Rg_\R(D)$ which we already proved, we conclude that $\Rg_\R(A) \simeq \Rg_\R(B)$.

       However we claim that there is no (max,+)-isomorphism from $\Rg(A)$ to $\Rg(B)$
      which sends the rays generated by the columns of $A$ to the rays generated by those of $B$, i.e.\ satisfying \Cref{thm:max_plus_isophi_general}~\ref{it_J_e_x}
      To show this, we shall use the fact that a (max,+)-isomorphism must be an isometry
      in Hilbert's seminorm.
      Recall that the Hilbert's seminorm on $\R^n$ is defined by 
      \[
      \|z\|_H = \max_i z_i -\min_j z_j  .
      \]
      or equivalently
      \[
      \|z\|_H =\inf \{\beta-\alpha \mid \alpha,\beta \in \R, \alpha e \leq z\leq \beta e \},
      \]
      where $e$ is the identically one vector.
      Using this property, it is clear that every (max,+)-isomorphism is nonexpansive in Hilbert's
      seminorm. Hence, a (max,+)-isomorphism between two subsets of $\R^n$ and $\R^p$ stable by the addition
      of constants must be an isometry in this seminorm, as claimed above.
Now, we check that
      $\|A_{\cdot 1}-A_{\cdot 2}\|_H = 4$,
      $\|A_{\cdot 1}-A_{\cdot 3}\|_H =\|A_{\cdot 2}-A_{\cdot 3}\|_H  =2$,
      while 
      $\|B_{\cdot 1}-B_{\cdot 2}\|_H = 4$,
      $\|B_{\cdot 1}-B_{\cdot 3}\|_H =3$ and $\|B_{\cdot 2}-B_{\cdot 3}\|_H  =1$,
      hence, the sets of columns of $A$ and $B$ cannot be put in isometric correspondence.

      This example is illustrated in~\Cref{fig-flag}. Observe that $\Rg_\R(A)$ is a subset of $\R^3$ which is
      invariant by the translation by the constant vector $(1,1,1)$. So to vizualize it, it suffices
      to draw the cross-section of $\Rg_\R(A)$ by the hyperplane orthogonal to the vector $(1,1,1)$.
This cross-section is shown on the figure.
      One can see that $\Rg_\R(A)$ and $\Rg_\R(B)$ have the same diameters
      in Hilbert's seminorm (the diameter is achieved by the distance between the two first columns, in each case). It is
      not difficult to check that the diameter is the only invariant for (max,+) spaces with two finite generators belonging to some $\R^n$,
      meaning that two such spaces are isomorphic iff they have the same diameter, which explains
      how this example was constructed. Note also that $A_{\cdot 1},A_{\cdot 2}$ are sup-irreducible and relatively-inf-irreducible, but that $A_{\cdot 3}=\inf^{\Rg(A)}((1,0,1),(0,1,1))=\sup(A_{\cdot 1},A_{\cdot 2})$ is neither, even though $\R e+\{A_{\cdot 1},A_{\cdot 2},A_{\cdot 3}\}$ are the only inf-irreducible elements in $\Rg(A)$. This underlines that the notions of irreducibility are distinct.

    \end{Remark}

        Nevertheless when the Dirac masses belong to the sets, every (max,+)-isomorphism is affine.
   
        \begin{Corollary}
        \label{thm:max_plus_isophi_general_ss_hyp_simple}
        Let $\G$ (resp.\ $\F$) be a complete subspace  of $\pmRX$ (resp.\ $\pmRY$) for some sets $\X$ and $\Y$, with both $\F$ and $\G$ satisfying \Cref{assumption:proper}. Let $J$ be a (max,+)-isomorphism from $\F$ onto $\G$. Assume that $\{\delta^\top_{y}\}_{y\in\Y}\subset \F$ and $\{\delta^\top_{x}\}_{x\in\X}\subset \G$. Then the statements of \Cref{thm:max_plus_isophi_general} hold, in particular there exists $g:\X\to\R$ and a bijective $\phi:\X\to\Y$ such that $Jf=g+f\circ \phi$.
        \end{Corollary}
          \begin{proof} Fix $x\in\X$ and $y\in\Y$, $\la\in\R$ and $\la'\in\R$. By definition \eqref{eq:def_e_x}, we have that $e_x=\delta^\top_{x}$ and $e'_y=\delta^\top_{y}$. Moreover $[\delta^\top_{x}+\la]$ is maximal since if there exists $h\in\G$ and $\alpha\in\R$ such that $\delta^\top_{x}+\la\le h+\alpha$, then $h=\delta^\top_{x}+h(x)$, similarly $[\delta^\top_{y}+\la']$ is maximal. Moreover \Cref{lem:prop_e_x} guarantees that these are the only maximal Archimedean classes, and \Cref{claim:archimd_preserved} shows that $J$ sends an $f\in\F$ with $[f]$ maximal onto a $Jf\in\G$ with $[Jf]$ maximal. Hence we ``send Dirac masses onto Dirac masses'' and \Cref{thm:max_plus_isophi_general}~\ref{it_J_e_x} holds.
          
          \end{proof}
          We will often meet in the sequel sets containing the Dirac masses. The other simple sets are the ones behaving like Lipschitz functions:
    \begin{Remark}[When $\G$ is inf and sup-stable]\label{rmk:inf-sup_stable_sets} We have shown in \cite[Proposition 6.2]{aubin2022tropical} that complete subspaces $\G\subset\R^\X$ that are both inf and sup-stable can always be described as the set of 1-Lipschitz functions for the weak $\pmR$-valued metric $\tilde d:\X\times \X \to \pmR$ with $\tilde d(z,x)=e_x(z)$. For inf-stable sets $\G$, the $e_x$ are (relatively)-inf-irreducible functions. It is possible sometimes to show they are the only ones even without resorting to Archimedean classes. For instance, for the set $\G_{\nearrow}^\R$ of increasing functions from $\pmR$ to $\pmR$, consider a candidate inf-irreducible $f$  and take $x_0\in\Dom(f)\cap \R$, $A_0=\{x\in\R\,|\, f(x)> f(x_0)\}$, then $f=\inf[\sup(f(x_0),f),f+\delta_{A_0}]$ where $\delta_{A_0}$ is the indicator function of $A_0$. This forces $f=e_{x'_0}+\la_0$ for some $x'_0\in\pmR$, and for $x'_0=+\infty$ we recover the inf-irreducible constant functions, which are Busemann points. So \Cref{thm:max_plus_isophi_general}-\ref{it_J_e_x} holds for all (max,+)-isomorphisms $J:\G_{\nearrow}^\R\to\G_{\nearrow}^\R$ since they preserve inf-irreducibility. This is a case where there are no maximal Archimedean classes as discussed in Remark~\ref{rmk:no_max_archimedean}. As we will see in \Cref{sec:Busemann_Lip} and hinted at already in \Cref{lem:prop_e_x}, maximal Archimedean classes when they exist allow to separate within inf-irreducible points the $e_x$ from the horoboundary elements, among which the Busemann points.%
    \end{Remark}

   \subsection{Sup-irreducible functions and a Krein--Milman theorem}\label{sec:c-convex}

    We have explored the inf-irreducible elements of a sup-stable $\G$ by studying the $e_x$, can we symmetrize the arguments to obtain the sup-irreducible elements? In general not, because $\G$ has no reason to be inf-stable. Nevertheless we can use the notion of minimal element of a subset $\A\subset \G$, i.e.\ the $u\in\A$ such that for any $v\in\A$ with $v\le u$, we have $v=u$. We say that $\G$ is \emph{stable by decreasing nets} if the infimum of any decreasing net of elements of $\G$ also belongs to $\G$. Considering minimal elements gives us a Krein--Milman theorem.
    \begin{Proposition}[A tropical Krein--Milman]\label{fct:extreme_pts_generators_minSx}%
		Let $\G$ be a complete subspace of $\pmRX$ and assume that $\G$ is stable by decreasing nets. For all $x\in \X$, define $S(x):=\{ u(\cdot)\in \G \, | \, u(x) \ge 0 \}$, and let $\Min(S(x))$ be the set of the minimal elements of $S(x)$. Then $\Min(S(x))$ is non-empty and, for any $u\in \Min(S(x))$ and $\la\in\R$, $u(\cdot)+\la$ is sup-irreducible. Moreover we have that these functions sup-generate $\G$, i.e.\ $\G=\{\sup_{x\in \X, u\in \Min(S(x))} u(\cdot) + a_{x,u}\,|\, a_{x,u}\in \mR \}$.
    \end{Proposition}

    \begin{proof} Since $\G$ is stable by decreasing nets, it follows
      from Zorn's lemma that $\Min(S(x))$ is non-empty. Take $u\in \Min(S(x))$. Notice that necessarily $u(x)=0$, as otherwise $(u(\cdot)-u(x))\in S(x)$ is smaller than $u$. Fix $x\in \X$, $f,g \in \G$ such that $u=\sup(f,g)$, then w.l.o.g.\ $f(x)=u(x)=0$, so $f\le u$ and $f\in S(x)$, hence $f=u$ and $u$ is sup-irreducible.
 
    Fix $f\in \G$. Since $\G$ is stable by decreasing nets, for every $x\in\X$ we can select $u\in \Min(S(x))$ such that for all $x'\in \X$, $f(x')-f(x)\ge u(x')$. In other words, for all $x'\in\X$,
    \begin{equation*}
        f(x') \ge \sup_{x\in \X, u\in \Min(S(x)), \, u\le f-f(x)} [u(x')+f(x)] \stackrel{x=x'}{\ge
        } f(x').
    \end{equation*}
    So  $\G=\{\sup_{x\in \X, u\in \Min(S(x))} u(\cdot) + a_{x,u}\,|\, a_{x,u}\in \mR \}$. 
	\end{proof}
    We say that $b:\X\times\X\to\pmR$ is strictly tropically monotone if $b(x,x)+b(y,y)\ge b(x,y)+b(y,x)$ for all $x,y\in\X$, with equality iff $x=y$. This holds for $(\X,d)$ a metric space with $b=-d$ or for $X$ a Hilbert space with $b=\bracket{\cdot,\cdot}$.
    
    \begin{Proposition}[A tropical inverse Krein--Milman]\label{lem:isophi_extreme_pts_generators_compact}
	Assume that $\Y$ is a compact Hausdorff space and $b(x,\cdot)$ is continuous over $\Y$ for all $x\in \X$. Then we have that every sup-irreducible element of $\Rg(B)$ belongs to $\{b(\cdot,y)+\lambda\}_{y\in\Y,\lambda\in\pmR}$. 
 
    If moreover $\Rg(B)$ is stable by decreasing nets, $\X=\Y$ and $b$ is strictly tropically monotone, then, $\{u(\cdot)+\la\}_{x\in X,\la \in \pmR, u\in \Min(S(x))}= \{b(\cdot,y)+\lambda\}_{y\in\Y,\lambda\in\pmR}$ and the $\{b(\cdot,y)+\lambda\}_{y\in\Y,\lambda\in\pmR}$ are the only sup-irreducible functions.
	\end{Proposition}
	\begin{proof}

        Take a sup-irreducible $f\in\Rg(B)$ with $f\notin \{x\mapsto+\infty,\, x\mapsto-\infty\}$. We know that $f=\bB \bB^\circ f$ \cite[p.3]{Akian04setcoverings} and want to show that $f=b(\cdot,y)+\la$ for some $y$ and $\la$. Define
        \begin{equation*}
             Z=\{\Z\subset \Y\,|\, \Z \text{ compact}, 
                f=\sup_{y\in\Z} b(\cdot,y)-\bB^\circ f(y)\}.
        \end{equation*}
        The set $Z$ is nonempty since $\Y\in Z$.
        We next show that it is inductive in the sense of Bourbaki, i.e.\ any totally ordered subset has a maximal element (here a minimal one).
        Consider a decreasing net $(\Z_\alpha)_{\alpha\in\A}\in Z^\Asc$.
        The intersection of compact sets $\bigcap_{\alpha\in\A} \Z_\alpha\in Z$ is compact. Moreover, by definition of $\bB^\circ$, the function $(\bB^\circ f=\sup_{x\in \X} b(x,\cdot)- f(x))$ is the supremum of continuous functions over $\Y$, so it is lower semicontinuous. Hence the function $(b(x,\cdot)- \bB^\circ f)$ is upper semicontinuous. By compactness
        \begin{equation*}
            \forall x\in\X, \, \forall \alpha\in\Asc, \, \exists y_{\alpha,x}\in\Z_\alpha, \, f(x)=b(x,y_{\alpha,x})-\bB^\circ f(y_{\alpha,x}).
        \end{equation*}
        Fix for each $x$, one such $y_{\alpha,x}$. Let $\Z_\cap:=\bigcap_{\alpha\in\A} \Z_\alpha$. Since the net is decreasing with respect to the inclusion, $\{y_{\alpha,x}\}_{\alpha\ge \alpha_0}\subset \Z_{\alpha_0}$, which is compact, so w.l.o.g.\ $y_{\alpha,x}\to y_x\in \Z_{\alpha_0}$. As this holds for any $\alpha_0$, by considering subnets, we obtain that $y_x\in \Z_{\cap}$. As $(b(x,\cdot)- \bB^\circ f)$ is upper semicontinuous,
        \begin{equation*}
            \forall x\in\X, \, f(x)=\lim_\alpha b(x,y_{\alpha,x})-\bB^\circ f(y_{\alpha,x}) \le b(x,y_{x})-\bB^\circ f(y_{x}).
        \end{equation*}
        However $f\ge \sup_{y\in\Z_\cap} b(\cdot,y)-\bB^\circ f(y)$. So $\Z_\cap\in Z$ and by Zorn's lemma, there exists a minimal $\bar \Z\in Z$. Assume that we can take $y_1,y_2\in\bar\Z$ with $y_1\neq y_2$. Then since $\Y$ is Hausdorff, we can find an open set $U\subset \Y$ such that $y_1\in U$ and $y_2\notin \bar U$. But this means that 
        \begin{equation*}
             f=\max\left(\sup_{y\in\bar\Z\cap \bar U} b(\cdot,y)-\bB^\circ f(y),\sup_{y\in \bar\Z\cap (\Y\backslash U)} b(\cdot,y)-\bB^\circ f(y)\right).
        \end{equation*}
        Since $f$ is sup-irreducible, it is equal to one of the two functions, but this contradicts the minimality of $\bar \Z$. Hence $\bar \Z$ contains a single element, and $f=b(\cdot,y)-\bB^\circ f(y)$ for $y\in \bar \Z$.
    
        Assume now that $\Rg(B)$ is stable by decreasing nets, $\X=\Y$ and $b$ is strictly tropically monotone. We showed in \Cref{fct:extreme_pts_generators_minSx} that the functions $u(\cdot)+\la$ were sup-irreducible, which gives one inclusion. To show the reverse inclusion, fix $x\in\X$, and take $u(\cdot)\in \Min(S(x))$ such that $b(\cdot,x)-b(x,x)\ge u(\cdot)$. Since $u\in\Rg(B)$, for any $\epsilon>0$, we can find $y_{x,\epsilon}$ such that $u(x)\ge b(x,y_{x,\epsilon})- (\bB^\circ u)(y_{x,\epsilon}) \ge  u(x)-\epsilon=-\epsilon$. So
        \begin{equation*}
            b(\cdot,x)- b(x,x) \ge u(\cdot) \ge b(\cdot,y_{x,\epsilon})- (\bB^\circ u)(y_{x,\epsilon}) \ge b(\cdot,y_{x,\epsilon})-b(x,y_{x,\epsilon}) -\epsilon.
        \end{equation*}
        Let $\epsilon$ go to $0$. By compactness of $\Y$, $y_{x,\epsilon}$ has a converging subnet converging to a $y_{x,0}$. Then, since $b(\textunderscore,\cdot)$ is continuous, evaluating the last expression at $y_{x,0} \in \X$, we obtain
        \begin{equation*}
            b(y_{x,0},x)- b(x,x) \ge b(y_{x,0},y_{x,0})-b(x,y_{x,0}), \quad \text{and} \quad u(\cdot)\ge b(\cdot,y_{x,0})-b(x,y_{x,0}).
        \end{equation*}
        The first part is precisely the reverse inequality of tropical monotonicity, so $y_{x,0}=x$. So $b(\cdot,x)- b(x,x)=u(\cdot)$ belongs to $\Min(S(x))$.
        
	\end{proof}
    When $\G$ is stable by decreasing nets, the first part of \Cref{lem:isophi_extreme_pts_generators_compact} shows that considering a subset $\bar \Y\subset \Y$ such that $\{u(\cdot)+\la\}_{x\in X,\la \in \pmR, u\in \Min(S(x))}= \{b(\cdot,y)+\lambda\}_{y\in\bar \Y,\lambda\in\pmR}$ would give a sup-irreducible family $\{b(\cdot,y)\}_{y\in\bar \Y}$ of $\Rg(B)$, which sup-generates it owing to \Cref{fct:extreme_pts_generators_minSx}. However we cannot ensure that $\bar \Y$ is compact, a problem that obviously does not happen in finite dimensions. Hence we will make sup-irreducibility an assumption. 
    
    We say that a kernel $b:\X\times \Y\to\R$ is fully-reduced if, for all $x\in\X$ and $y\in\Y$, $b(x,\cdot)$ and $b(\cdot,y)$ are sup-irreducible and if, for all $x_0,x_1\in\X$, $y_0,y_1\in\Y$ and $\la\in\R$, the implications $(b(\cdot,y_0)=b(\cdot,y_1)+\la \implies y_0=y_1)$ and $(b(x_0,\cdot)=b(x_1,\cdot)+\la \implies x_0=x_1)$ both hold, i.e.\ the rows and columns of $B$ are pairwise independent. It is sometimes possible to show this property directly as in \Cref{fct:cvx_extreme_pts} for convex lower semicontinuous functions. Otherwise the second part of \Cref{lem:isophi_extreme_pts_generators_compact} provides a testable inequality, satisfied for instance for any $p>0$ by $b=-d(x,y)^p$ over a metric space $(\X,d)$. For $p\in(0,1]$, $b$ defines the inf-stable $\Rg(B)$ of Hölder functions, which are indeed stable by decreasing nets. If $\Y$ is finite, $\Rg(B)$ is also stable by decreasing nets.%

    \begin{Theorem}\label{thm:rg_B_isometry} Let $\X,\X',\Y,\Y'$ be compact Hausdorff topological spaces, and $b:\X\times \Y\to\R$, $c:\X'\times \Y'\to\R$ be two continuous fully-reduced kernels. Then the following statements are equivalent:
		\begin{enumerate}[label=\roman*),labelindent=0cm,leftmargin=*,topsep=0.1cm,partopsep=0cm,parsep=0.1cm,itemsep=0.1cm]
			\item \label{it_rgB_isophi} there exists a (max,+)-isomorphism $J:\Rg(B)\to\Rg(C)$;
            \item \label{it_rgB_circ_isophi} there exists a (max,+)-isomorphism $J^\circ:\Rg(B^\circ)\to\Rg(C^\circ)$;
			\item \label{it_kernel_isometry} the two kernels are such that there exist two homeomorphisms $\tau:\X'\to\X$ and $\sigma:\Y'\to\Y$, and two continuous functions $\psi:\X'\to\R$ and $\varphi:\Y'\to\R$ satisfying
       \begin{equation}\label{eq:kernel_isometry}
           c(x',y')=\psi(x')+b(\tau(x'),\sigma(y'))+\varphi(y').
       \end{equation}
       \item \label{it_rgB_isophi_any} all (max,+)-isomorphisms $J:\Rg(B)\to\Rg(C)$ are of the form $Jf=\psi+f\circ \tau$ for some $\psi:\X'\to\R$ and a bijective $\tau:\X'\to\X$, and there exists one such $J$;%
       \item \label{it_rgB_circ_isophi_any} all (max,+)-isomorphisms $J^\circ:\Rg(B^\circ)\to\Rg(C^\circ)$ are of the form $J^\circ h=\varphi+h\circ \sigma$ for some $\varphi:\Y'\to\R$ and a bijective $\sigma:\Y'\to\Y$, and there exists one such $J^\circ$.
		\end{enumerate}
	\end{Theorem}%
    \begin{proof}
        \ref{it_rgB_isophi}$\Leftrightarrow$\ref{it_rgB_circ_isophi} Given $J$, set $J^\circ=\bar C J \bB$. Conversely, given $J^\circ$, set $J=\bar C^\circ J^\circ \bB^\circ$.
        
        \ref{it_rgB_isophi}$\Rightarrow$ \ref{it_rgB_isophi_any} $\Rightarrow$\ref{it_kernel_isometry} Let $y'\in\Y'$. Since $c(\cdot,y')$ is assumed to be sup-irreducible, its image by $J^{-1}$ is sup-irreducible, and thus necessarily of the form $b(\cdot,\sigma(y'))+\varphi(y')$ by \Cref{lem:isophi_extreme_pts_generators_compact}. Since the kernel functions are the only sup-irreducible elements, their images by anti-isomorphisms are the only relatively-inf-irreducible elements of the considered sets. These images are precisely the $e_x$ by Remark~\ref{rmk:e_x_et_b_x}, so $J$ maps $e_x$ on $e'_{x'}$. We can thus use \Cref{thm:max_plus_isophi_general}, and obtain that $Jf=\psi+f\circ \tau$ for some $\psi:\X'\to\R$ and a bijective $\tau:\X'\to\X$, whence \ref{it_rgB_isophi_any} holds. Moreover%
        \begin{equation*}        c(\cdot,y')=J(b(\cdot,\sigma(y'))+\varphi(y'))=\psi(\cdot)+b(\tau(\cdot),\sigma(y'))+\varphi(y').
        \end{equation*}
        Let us show the regularity of the functions, but let us avoid dealing with all functions simultaneously. We will use a translation property. Fix a $y'_0\in\Y'$. As in \Cref{lem:extension_order}, $J$ allows to define a (max,+)-isomorphism $J'$ from $\Rg(B')$ to $\Rg(C')$ where $b'(x,\sigma(y')):=b(x,\sigma(y'))-b(x,\sigma(y_0'))$ and $c'(x,y'):=c(x,y')-c(x,y_0')$. Since $J'$ also satisfies \Cref{thm:max_plus_isophi_general}-\ref{it_J_e_x}, we can thus use \eqref{eq:charac_J}, and obtain that $J'f=\psi'+f\circ \tau$ with the same $\tau$ and $\sigma$ as $J$ and some $\psi'(\cdot)$, owing to \Cref{lem:extension_order}. Moreover, as previously because sup-irreducibility is preserved by translations, we must have $c'(\cdot,y')=J'(b'(\cdot,\sigma(y'))+\varphi'(y'))$. Hence
        \begin{equation*}
        c'(\cdot,y')=J'(b'(\cdot,\sigma(y'))+\varphi'(y'))=\psi'(\cdot)+b'(\tau(\cdot),\sigma(y'))+\varphi'(y').
        \end{equation*}
        In particular, evaluating at $y'=y'_0$, by definition of $b'$ and $c'$, we have that $\psi'(\cdot)=-\varphi'(y'_0)$, so $c'(\cdot,y'_0)=b'(\tau(\cdot),\sigma(y'_0))$. Let $(x'_\alpha)_{\alpha\in\Asc}$ be a net converging to some $x'\in\X'$. Assume that $(\tau(x'_\alpha))_{\alpha\in\Asc
        }$, which belongs to the compact set $\X$, has two subnets converging respectively to $x_0$ and $x_1\in\X$. Then taking the limit, since the kernels are continuous we obtain that
        \begin{equation*}
        b'(x_0,\sigma(y'))+\varphi'(y')-\varphi'(y_0')=c'(x',y')=b'(x_1,\sigma(y'))+\varphi'(y')-\varphi'(y_0'),
        \end{equation*}
        so $b'(x_0,\sigma(y'))=b'(x_1,\sigma(y'))$. Hence by definition of $b'$, for all $y'\in\Y'$,
         \begin{equation*}
        b(x_0,\sigma(y'))-b(x_0,\sigma(y_0'))=b(x_1,\sigma(y'))-b(x_1,\sigma(y_0')).
        \end{equation*}
        However since the kernel $b$ is fully-reduced and $\sigma$ is bijective, we must have $x_0=x_1$. Consequently, all converging subnets have the same limit $x_0$. Similarly to \cite[Lemma 2.17]{aliprantis2006}, if $(\tau(x'_\alpha))_{\alpha\in\Asc
        }$ did not converge to $x_0$, then it would have a subnet staying away from a neighborhood of $x_0$, but by compactness this subnet has itself a converging subnet to $x_0$. We have thus shown that $\tau$ is continuous. Symmetrizing the proof, so is $\tau^{-1}$ and we get that $\tau$ and $\sigma$ are homeomorphisms. As we already established \eqref{eq:kernel_isometry}, since $b$ and $c$ are continuous and the bijections $\sigma$ and $\tau$ are continuous, so are $\psi$ and $\varphi$.

        \ref{it_rgB_circ_isophi}$\Rightarrow$ \ref{it_rgB_circ_isophi_any} $\Rightarrow$\ref{it_kernel_isometry} This implication is shown similarly, symmetrizing the arguments.
        
        \ref{it_rgB_isophi}$\Leftarrow$\ref{it_kernel_isometry} Just set $Jf=\psi+f\circ \tau$ defined for $f\in \Rg(B)$. Since $J:\Rg(B)\to\Rg(C)$ is clearly a (max,+)-isomorphism, it commutes with the operations in $f=\sup_{y\in\Y}b(\cdot,y)-\bB^\circ (y)$. From \eqref{eq:kernel_isometry} we deduce that $Jf\in\Rg(C)$, so we made one (max,+)-isomorphism explicit.
    \end{proof}

    \section{Studying (max,+) and order isomorphisms over specific sets}\label{sec:examples}
    
    Applying the results of \Cref{sec:general_sets}, and in particular \Cref{thm:max_plus_isophi_general}, we now characterize the (max,+)-isomorphisms for some given sets of functions: i) lower semicontinuous or continuous, ii) 1-Lipschitz with a (weak) locally compact or complete metric, iii) lower semicontinuous and convex. For some of these sets we can even characterize the order isomorphisms in great generality.
  
    \subsection{Indicator functions and lower semicontinuous functions}\label{sec:examples_continuous}
    
    The Banach-Stone theorem is a classical result of linear functional analysis. In this section, we give a (max,+) version of it. The following characterization \eqref{eq:isophi_lsc} appeared in \cite[Theorem 2.2]{KM} for (max,+)-isomorphisms over sets of continuous functions on separable locally compact topological spaces with the functions diverging to $\infty$ outside of a compact set, see \cite[Section 1.4, p34]{KM} for their notations. Inf-irreducible points were not used in the proof of Maslov and Kolokoltsov, the latter relied instead on a representation theorem of integral operators which we circumvented.%
    \begin{Theorem}[(max,+)-isomorphisms of lsc]\label{thm:max_plus_isophi_lsc} 
		Let $\G$ (resp.\ $\F$) be the space of lower semicontinuous functions over a Hausdorff topological space $\X$ (resp.\ $\Y$). Then every (max,+)-isomorphism $J$ from $\F$ onto $\G$ is of the form
		\begin{equation}\label{eq:isophi_lsc}
			Jf(x)=g(x)+f(\phi(x))
		\end{equation}
		where $g:\X \rightarrow \R$ is a continuous function and $\phi:\X\rightarrow\Y$ is a homeomorphism. For $\X$ and $\Y$ compact, the same holds if lower semicontinuous is replaced by continuous, or if the functions are restricted to be proper.%
        
	\end{Theorem}
    \begin{proof}
        By \Cref{thm:max_plus_isophi_general_ss_hyp_simple},  since $e_x=\delta_x^\top$, we already have that \eqref{eq:isophi_lsc} holds for some $g\in\R^\X$ and a bijective $\phi$. We just have to characterize more finely the properties of $g$ and $\phi$.  We have that $J(x\mapsto 0)=g$ and $J^{-1}(x\mapsto 0)=-g\circ \phi^{-1}$. Since $J$ and $J^{-1}$ are bijective over lower semicontinuous functions, $g$ and $-g\circ \phi^{-1}$ are lower semicontinuous. We use now that lower semicontinuous functions are stable by addition to show that $\phi$ is a homeomorphism. Let $F$ be a closed subset of $\Y$, then $\delta_F$ is lower semicontinuous, and so is $J(-g\circ \phi^{-1}+\delta_F)\stackrel{\eqref{eq:isophi_lsc}}{=}\delta_{\phi^{-1}(F)}$ since $J$ is an isomorphism. So $\phi^{-1}(F)$ is thus closed in $\X$. Consequently $\phi$ is bijective and continuous. Repeating the proof, the same holds for $\phi^{-1}$, thus $\phi$ is a homeomorphism. As the image of the closed epigraph of $-g\circ \phi^{-1}$ by the homeomorphism $(\phi,\Id_\R)$ is closed, composing $-g\circ \phi^{-1}$ with $\phi$ gives that $-g$ is also lower semicontinuous, entailing that $g$ is continuous. Conversely any $J$ as in \eqref{eq:isophi_lsc} is a (max,+)-isomorphism from $\F$ to $\G$.

        Concerning the set of continuous functions $\C(\X,\R)$, for $\X$ compact, by \Cref{lem:extension_order}, since its sup-closure is the set of proper lower semicontinuous functions, the (max,+)-isomorphisms between spaces of continuous functions are of the form \eqref{eq:isophi_lsc}. The converse part is immediate.%
    \end{proof}

    \begin{Theorem}[order isomorphisms of lsc]\label{thm:max_isophi_lsc}
    Let $\G$ (resp.\ $\F$) be the space of lower semicontinuous functions over a Hausdorff topological space $\X$ (resp.\ $\Y$). Then every order isomorphism $J$ from $\F$ onto $\G$ is of the form
		\begin{equation}\label{eq:isophi_lsc2}
			Jf(x)=g(x,f(\phi(x)))
		\end{equation}
		where $\phi$ is a homeomorphism and $g:\X \times \pmR \rightarrow \pmR$ is jointly lower semicontinuous, and $g(x,\cdot)$ bijective and increasing for all $x\in\X$ with inverse $g^1(x,\cdot)$ such that $g^1(\cdot,\cdot)$ is also jointly lower semicontinuous. For $\X$ and $\Y$ compact, the same holds if lower semicontinuous is replaced by continuous and $J$ satisfies \eqref{eq:J_continuity}.

	\end{Theorem} 
    
    Concerning lattice isomorphisms $J$ of continuous functions over compact Hausdorff spaces, assuming $J$ to be continuous under pointwise convergence, which is weaker than \eqref{eq:J_continuity}, formula \eqref{eq:isophi_lsc2} was first shown in \cite[Theorem 1]{Kaplansky1948} with an elegant proof based on characterizing ``prime ideals'' of $C(\X,\R)$. When the continuity assumption is dropped, formula \eqref{eq:isophi_lsc2} still holds on a $G_\delta$-dense set \cite[Theorem 1, Lemma 4]{Sanchez2008}, which encompasses the results of Kaplansky. In  \cite[Theorem 4.5]{Leung2016}, the same expression was obtained for isomorphisms of some specific vector subspaces of continuous functions on a metric space, named ``near vector lattices'' by the authors, and encompassing (locally) Lipschitz or uniformly continuous (bounded) functions. Compactifications of $\X$ and continuity of $J$, which is not automatic, are key elements in these approaches. They do not appear when focusing on lower semicontinuous functions.
    \begin{proof} The $\G$-relatively-inf-irreducible points of lower semicontinuous functions are the Dirac masses $(\delta_x+\la)_{x\in\X,\la\in\pmR}$. Indeed if $g,h\in\pmRY$ are such that $\delta_x=\infG(g,h)$, then $\delta_x\le g$ and $\delta_x\le h$, so both are Dirac masses $\delta_x+\la$, one being equal to $\delta_x$. Conversely take a $\G$-relatively-inf-irreducible $f\in\G$. It is necessarily inf-irreducible by \Cref{claim:inf-relatif}. Assume there exists $x_0\neq x_1$ such that $f(x_0)< \infty$ and $f(x_1)< \infty$, then fix a closed subset $F\subset \X$ such that $x_0\in F$, $x_0\notin \overline{\X\backslash F}$ and  $x_1\notin F$. Hence $f=\inf(f+\delta_F,f+\delta_{\overline{\X\backslash F}})$, both of these functions being different than $f$. Thus $f=\delta_x+\la$ for some $\la\in\pmR$.
    
    Fix $y\in\Y$ and $\la \in \R$. Using \Cref{lem:isophi_extreme_pts}, we have that $J$ maps $\delta_y+\la$ on some $\delta_{\phi(y,\la)}+g_y(\la)$, for some $\phi:\Y \times \R \rightarrow\X$ and $g:\Y \times \R \rightarrow\R$, the $\infty$-valued functions being mapped onto themselves.

        Take $\la,\la'\in\R$, w.l.o.g.\ we can assume that $\la\le \la'$, so $J(\delta_y+\la)\le J(\delta_y+\la')$ which means that the two Dirac functions are comparable, hence $\phi(y,\la)=\phi(y,\la')$, and $g_y(\la)\le g_y(\la')$. As $J$ is bijective, the increasing $g_y:\R\rightarrow\R$ has to be bijective too, whence it is continuous. We set $\phi(y)=\phi(y,\cdot)$ from now on. Let $y_0\in\Y$ be such that $J\delta_{y_0}=\delta_{\phi(y)}+\la'$ for some $\la'\in\R$. Then we can find $\la\in\R$ such that $\la'=g_y(\la)$, so $J\delta_{y_0}=J(\delta_y+\la)$, hence $y_0=y$ and $\la=0$, in other words $\phi:\Y \rightarrow \X$ is bijective. Similarly we can show that
        \begin{equation*}
            J^{-1}(\delta_x+\la)=\delta_{\phi^{-1}(x)}+g^1_{\phi^{-1}(x)}(\la)
        \end{equation*}
        for some increasing function $g^1_y:\R\rightarrow\R$ satisfying $g^1_y(g_y(\la))=\la$, since $J^{-1}J=\Id_\F$.

        Take any l.s.c.\ $f\in\pmRY$ and $h\in\pmRX$, fix $y\in\Y$ and $x\in\X$. Then $f \le \delta_y +f(y)$ and $h \le \delta_x +h(x)$. So composing with $J$ and $J^{-1}$, the following holds
        \begin{align*}
            (Jf)(\phi(y)) &\le g_y(f(y)),\\
            (J^{-1}h)(\phi^{-1}(x)) &\le g^1_{\phi^{-1}(x)}(h(x)).
        \end{align*}
        Replace $h$ by $Jf$ and $\phi^{-1}(x)$ by $y$ in the second inequality, thus
        \begin{equation*}
            f(y)\le g^1_{y}(Jf(\phi(y)))\le g^1_{y}(g_y(f(y)))=f(y).
        \end{equation*}
        Composing with $g_y$, we obtain $Jf(\phi(y))=g_y(f(y))$. Renaming variables to write $Jf(x)=g_x(f(\phi(x)))$ and $J^{-1}f(x)=g^1_x(f(\phi^{-1}(x)))$, we get the desired result \eqref{eq:isophi_lsc2} with now $g_x(g^1_{\phi(x)}(f(x)))=f(x)$ after renaming. Moreover $J(y\mapsto \la)=g(\cdot,\la)$ is lower semicontinuous. The same holds for $g^1$. Thus $g(\cdot,\cdot)$ is separately lower semicontinuous with $g(x,\cdot)$ increasing. Consequently, by \cite[Theorem 1]{Grushka2019}, $g(\cdot,\cdot)$ is jointly lower semicontinuous.\footnote{ \cite[Theorem 1]{Grushka2019} is written for continuity rather than lower semicontinuity, but taking $T=T_1=\R$ in Grushka's notations, and setting $\tau_2=+\infty$ in his proof, the result straightforwardly extends to lower semicontinuous functions since they have almost the same characterization through neighborhoods as continuous functions.} By similar arguments, so is $g^1$. Using only the indicator functions, we have obtained the intermediary result that
        \begin{equation}\label{eq:isophi_lsc_inter}
            Jf(x)=g(x,f(\phi(x)))
        \end{equation}
        with $\phi$ bijective and $g$ jointly lower semicontinuous, with $g$ continuous and increasing in the second variable.
        
        Let $F$ be a closed subset of $\Y$, then $\delta_F$ is lower semicontinuous, so the sum $g^1(\cdot,0)+\delta_F(\cdot)$ is also lower semicontinuous. We have that $J(g^1(\cdot,0)+\delta_F(\cdot))(x)=g_x(g^1(\phi(x),0))=0$ for $x$ such that $\phi(x)\in F$ and $\infty$ otherwise, since $g_x(\infty)=\infty$ by increasingness and bijectivity of $g_x(\cdot)$. Hence $J(g^1(\cdot,0)+\delta_F(\cdot))=\delta_{\phi^{-1}(F)}$ is lower semicontinuous and $\phi^{-1}(F)$ is closed. Consequently $\phi$ is bijective and continuous, repeating the proof, the same holds for $\phi^{-1}$, thus $\phi$ is a homeomorphism.

        The converse part of a theorem is just a consequence of the following claim.
        \begin{Claim}
            Let $g:\X \times \R \rightarrow \R$ be jointly lower semicontinuous and bijective and strictly increasing in the second argument. Take any lower semicontinuous $f$. Then $g(\cdot, f(\cdot))$ is lower semicontinuous.
        \end{Claim}
        \begin{proof}
             Fix a sequence $(x_n)_{n\in\N}$, converging to some $\xb$,
             \begin{align*}
                 g(\xb,f(\xb))\le \liminf_{n\rightarrow \infty} g(x_n,f(\xb))&\le  \liminf_{n\rightarrow \infty} g(x_n,\liminf_{m\rightarrow \infty}f(x_m))\\
                 &=\liminf_{n\rightarrow \infty} \liminf_{m\rightarrow \infty} g(x_n,f(x_m))\le \liminf_{n\rightarrow \infty} g(x_n,f(x_n)).
             \end{align*}
         where the first inequality stems from the lower semicontinuity of $g(\cdot,f(\xb))$, the second from that of $f$ and the monotonicity of $g(x_n,\cdot)$, the equality is a consequence of the continuity of $g(x_n,\cdot)$ which follows from the bijectivity and increasingness of $g(x_n,\cdot)$. The last inequality is a diagonal argument over the joint lower semicontinuity.
        \end{proof}

          Concerning the set of continuous functions $\C(\X,\R)$, for $\X$ compact and $J$ satisfying \eqref{eq:J_continuity}, by \Cref{lem:extension_order}, since its sup-closure is the set of proper lower semicontinuous functions, the order isomorphisms between spaces of continuous functions are of the form \eqref{eq:isophi_lsc_inter}. In all the proof above, we just replace lower semicontinuous by continuous with no other change to derive the result applying \Cref{lem:extension_order}.%
    \end{proof}

    About anti-involutions, if $\X$ is discrete, then $-\Id$ is an anti-involution over (non proper) lower semicontinuous functions. However if the space has no isolated points then there are no anti-involutions.
    
	\begin{Lemma}\label{lem:lsc_no_anti}
		Let $\G$ (resp. $\F$) be the space of lower semicontinuous functions over a set $\X$ (resp. $\Y$). If there exists an order anti-isomorphism $T$ from $\F$ to $\G$, then $\X$, the set $I_{\X}$ of the isolated points of $\X$, $\Y$ and $I_{\Y}$ are all in bijection. In particular, if $\X$ or $\Y$ do not have isolated points, or if $I_{\X}$ and $I_{\Y}$ do not have the same cardinality, then there is no order anti-isomorphism. %
	\end{Lemma}
    The intuition here is that anti-isomorphisms send relative-inf-irreducible points onto sup-irreducible points by \Cref{lem:isophi_extreme_pts}. However the sup-irreducible points of lower semicontinuous functions are only the $\delta^\bot_x$ of isolated $x$.
	\begin{proof}
        Let us first characterize the sup-irreducible points. Let $f$ be a lower semicontinuous function taking finite values at two distinct points $x_0$ and $x_1$. Define two functions $\chi_0,\chi_1$ by setting $\chi_0(x_0)=-1$, $\chi_1(x_1)=-1$ and $0$ elsewhere. Then $f=\sup(f+\chi_0, f+\chi_1)$, with the two latter functions being lower semicontinuous as a sum of such functions, so $f$ is not sup-irreducible. The same applies if $f$ takes the $+\infty$ value at two distinct points $x_0$ and $x_1$ by setting $f_0=f$ everywhere except at $x_0$ where we choose $f_0(x_0)=0$ (resp.\ $f_1(x_1)=0$), so $f=\sup(f_0,f_1)$. The function $f_0$ is l.s.c.\ since its epigraph is the union of that of $f$ and of $\{x_0\}\times[0,+\infty]$. This shows that every sup-irreducible function has to be of the form $\delta^\bot_x+\la$ with $\la\in \pmR$. For $\delta^\bot_x$ to be lower semicontinuous, we need to have $\X \backslash \{x\}$ to be closed, i.e.\ $x$ to be isolated.

        Assume there exists an anti-isomorphism $T$.  Similarly to the proof of \Cref{thm:max_isophi_lsc}, applying \Cref{lem:isophi_extreme_pts}, it is easy to show that $T(\delta^\top_y+\la)=\delta^\bot_{\phi(y)}+g(y,\la)$ with $\phi:\Y\to I_{\X}$ bijective and $g(y,\la)\in\R$. In the same fashion, we have a bijective $\psi:\X\to I_{\Y}$. Hence $\psi\circ\phi:\Y\to I_{\Y}$ is injective, and by the Cantor--Schröder--Bernstein theorem, $\Y$ and $I_{\Y}$ are thus in bijection. Considering $T^{-1}$, the sets $\X$ and $I_{\X}$ are also in bijection. So we can construct with these bijections and the previous injections, two injections from $\X$ to $\Y$ and $\Y$ to $\X$. Again by the Cantor--Schröder--Bernstein theorem, we derive a bijection between $\X$ and $\Y$.%

	\end{proof}
    Note that \Cref{lem:lsc_no_anti} does not rule out that there can be interesting anti-isomorphisms $T:\F\to \G$ over lower semicontinuous functions $\F$ with $\G$ another function space, as given for instance by \cite[Theorem 2.2]{Bachir2001} for $\X$ a Banach space. %

            \subsection{Busemann points and 1-Lipschitz functions}\label{sec:Busemann_Lip}%
   In \cite[Theorem 1]{Sanchez2011}, it was first shown that, for Lipschitz functions on bounded metric spaces, the structure of order isomorphisms was necessarily that of lower semicontinuous functions, obtained in \Cref{thm:max_isophi_lsc}, with a bi-Lipschitz $\phi$. Extending this result to general metric spaces, \cite[Theorem 5.5, eq.(8)]{Leung2016} provided an ``if an only if'' condition on the existence of a bi-Lipschitz $\phi$. However at the moment, to the best of our knowledge there is no sufficient condition on the form of the function $g:\X\times\R\rightarrow\R$, not even lower semicontinuity,\footnote{On Lipschitz functions, even for $\phi=\Id$, the monograph \cite[Chapter 7]{Appell1990} says: \emph{``Quite amazing is the fact that the operator $J$ may act in a Hölder space, although the generating function $g$ is not continuous.''} Further conditions, like global Lipschitzianity of $J$, can however ensure that $J$ is affine as in \cite[Theorem 7.8]{Appell1990} for $\X=[0,1]$. The necessity of extra assumptions on $J$ is also discussed therein if one is interested by differentiability properties.} nor any specific results for order isomorphisms of the space of $1$-Lipschitz functions.

    Nevertheless a similar result to what follows was shown in \cite[Corollary 4.1]{Cabello2017} for 1-Lipschitz functions over complete metric spaces, for order isomorphisms that also preserve convex combinations, i.e.\ $J(sf+(1-s)f')=sJf+(1-s)Jf'$, through very different proof techniques such as bases of topologies. The line of research of Cabello et al.\ on weak metrics was adapted later in \cite[Theorem 3.1]{Daniilidis2020} to Finsler manifolds. We instead do not ask that convexity is preserved but that $J$ commutes with the addition of constants. 
            \begin{Definition}
            If $\X$ is a set, we say that a map $\delta: \X \times \X \to \R$ is a {\em weak metric}
              if,
               for all $x,y,z\in \X$, $\delta(x,z)\leq \delta(x,y)+\delta(y,z)$,
            and if the symmetrized
        map $\delta_s(x,y)=\delta(x,y)+\delta(y,x)$ defines an ordinary metric on $\X$. We say that a map $f: \X\to \R$ is {\em nonexpansive}
with respect to $\delta$, or $1$-Lipschitz, if $f(x)  \leq \delta(x,y)+ f(y) $ holds for all $x,y\in \X$.
We denote by $\Lipb(\X,\R)$ the set of $1$-Lipschitz maps $f:\X \to \R$.
            \end{Definition}

Unlike \cite{Cabello2017,Daniilidis2020}, we do not assume the weak metric to be nonnegative-valued, which allows us to consider the Funk metric in \Cref{example:Funk}, for which $\delta_s$ is the Hilbert metric. We still assume that $\delta_s$ is a metric, hence has nonnegative values. Our result reads as follows

\begin{Theorem}[(max,+)-isomorphisms of 1-Lip]\label{th-Lip}
Let $\G$ (resp.\ $\F$) be the set of nonexpansive maps over a set equipped with a weak metric $(\X,\delta)$ (resp. $(\Y,\delta')$). Assume that $(\X,\delta)$ satisfies either i) that the closed balls of $(\X,\delta_s)$ are compact or ii) that $\delta=\frac12 \delta_s$ and $(\X,\delta_s)$ is complete, and that the same is true for $(\Y,\delta')$. Then every (max,+)-isomorphism $J$ from $\F$ onto $\G$ is of the form
    \begin{equation}\label{eq:max_plus_isophi_Lip}
        Jf(x)=g(x)+f(\phi(x))
    \end{equation}
  where $g:\X\to\R$ is nonexpansive and $\phi:\X\rightarrow\Y$ is such that $g(x)-g(x')+\delta'(\phi(x),\phi(x'))=\delta(x,x')$ for all $x,x'\in\X$. If both $\delta$ and $\delta'$ are metrics, then $g$ is constant and $\delta'(\phi(x),\phi(x'))=\delta(x,x')$.
\end{Theorem}
For $p\in(0,1]$ and $\delta(x,y): \X \times \X \to [0,+\infty)$, by subadditivity of $t\in\R_+\mapsto t^p$, the result above also applies to $p$-Hölder functions since $\delta^p$ is also a weak metric. Other interesting examples of weak metrics are those of the form ``$ \inf_\gamma \int L(\gamma)$'' with nonnegative Lagrangian $L$, used in \cite[Chapter 7]{Villani2009} or \cite[Section 7.2]{aubin2022tropical}.

The proof of \Cref{th-Lip} uses the compactification of $(\X,\delta)$ by horofunctions and is delayed to p.\pageref{proof:th-lip}, after showing some intermediary lemmas. The horofunction compactification was
introduced by Gromov~\cite{gromov78}, see also~\cite{gromov}, \cite[Ch.~II]{ballmann} and~\cite{rieffel}.
It carries over to weak metrics, see~\cite{AGW,GV10,Walsh2018}. For compact $\X$, \Cref{th-Lip} is merely a specialization of \Cref{thm:rg_B_isometry}.
       
	\begin{Corollary}
	    Let $\G$ (resp.\ $\F$) be the set of $1$-Lipschitz functions over a complete metric space $(\X,d)$ (resp.\ $(\Y,d')$). Then every (max,+)-anti-isomorphism $J$ from $\F$ onto $\G$ is of the form
		\begin{equation}\label{eq:anti-inv_isophi_Lip}
			Jf(x)=\alpha-f(\phi(x))
		\end{equation}
		where $\alpha\in \R$ is a constant and $\phi:\X\rightarrow\Y$ is an isometry. 
	\end{Corollary}
 \begin{proof}
    Notice that $-\Id$ is an anti-involution over $1$-Lipschitz functions, and just apply \Cref{prop:isophi_anti-invol} in conjunction with \eqref{eq:max_plus_isophi_Lip}.
 \end{proof}

       Fix a {\em basepoint} $\bar x\in \X$.
        We define a map $\imath:\X\to\R^\X$ by associating
        to any $x\in \X$ the function $\imath(x) \in \Lipb(\X,\R)$,
        \[
        \imath(x)(y)= \delta(y,x)-\delta(\bar x,x)  .
        \]%
        By the triangle inequality, $-\delta(y,\bar x) \leq \imath(x)(y) \leq \delta(y,\bar x)$, so that
        $i(\X)$ can be identified to a subset of the Cartesian product $\prod_{y\in \X}[-\delta(y,\bar x), \delta(y,\bar x)]$
        which is compact in the product topology, by Tychonoff's theorem. The {\em horocompactification} of $\X$ is defined as
        the closure $\overline{\imath(\X)}$ of $\imath(\X)$  in the product topology. As a special case of a ``second'' Ascoli's theorem \cite[p310]{Schwartz1970-ew}, since $\Lipb(\X,\R)$ with the pointwise convergence topology is isomorphic to itself with the compact-open topology and thus with uniform convergence on compact sets, the closure of $\imath(\X)$
        is the same if we take the topology of uniform convergence on compact sets of the metric space $(\X,\delta_s)$.
        We set $\X(\infty)\coloneqq \overline{\imath(\X)}\setminus \imath (\X)$. An element of $\X(\infty)$ is called a {\em horofunction}.%

        Following \cite[Definition 2.1]{Walsh2018}, we say that a net $\{z_\alpha\}_{\alpha\in\Asc}\subset \R$ is \emph{almost non-decreasing} if, for any $\epsilon>0$, there exists $A$ such that $z_\alpha\le z_{\alpha'}+\epsilon$ for all $\alpha,\alpha'\in\Asc$ such that $A\le\alpha\le\alpha'$. As easily seen, taking limits, we have $\limsup_\alpha z_\alpha\le \liminf_\alpha z_{\alpha}-\epsilon$ for all $\epsilon>0$, so bounded almost non-decreasing nets converge in $\R$.%

        We say that a net $(x_\alpha)_{\alpha}$ of elements of $\X$ is an {\em almost-geodesic}~\cite{AGW,Walsh2018} if for all $\epsilon>0$,
        \[ \delta(\bar x, x_\beta) \geq \delta(\bar x, x_\alpha) + \delta(x_\alpha,x_\beta) - \epsilon
        \]
        holds for all $ \alpha $ and $\beta$ large enough, with $\alpha\leq \beta$.
        This definition is similar to the one of Rieffel~\cite{rieffel}, except that Rieffel requires
        an almost-geodesic to be a curve, see~\cite[Prop.~7.12]{AGW} for a comparison.
        Note that a restricted definition of almost-geodesics, defined as sequences, was adopted
        in~\cite{AGW}, the more general notion in terms of nets was introduced in~\cite[Definition 2.2]{Walsh2018}.
        When the closed balls of $(\X,\delta_s)$ are compact, then, $\X$ is a countable union of compact sets, and so, the topology of $\overline{\imath(\X)}$ is metrizable. Then, limits in this space can be characterized in terms of sequences, rather than of nets, and so it suffices to consider sequences that are
        almost-geodesic.

        We say that an horofunction $h\in \X(\infty)$ is a {\em Busemann point} iff it is the pointwise limit
        of a net $\imath(x_\alpha)$ where $x_\alpha$ is an almost-geodesic.
        We recall the following result, which was proved in~\cite{AGW}.
        See also~\cite{Walsh2018}.
        \begin{Theorem}[{\cite[Coro.~6.6 and 7.7]{AGW}}]\label{th-agw}
          The set of inf-irreducible elements of $\Lipb(\X)$ consists of the elements
          of $\imath(\X)$ and of the Busemann points of $\X(\infty)$, up to translations by additive constants.
        \end{Theorem}
        We shall abuse terminology and say that an almost-geodesic $x_\alpha$ {\em converges} to a horofunction $h$ to mean that $\imath(x_\alpha)$ tends to $h$.
        \begin{Proposition}
          \label{prop-infty}%
Assume either that i) all the closed balls of $(\X,\delta_s)$ are compact or ii) $\delta=\frac12 \delta_s$ and $(\X,\delta_s)$ is complete. Then, if $h$ is a Busemann point,
then, any almost-geodesic $x_\alpha$ converging to $h$ is such that $\delta_s(\bar x, x_\alpha)$
tends to infinity.
        \end{Proposition}
        As already discussed in Remark~\ref{rmk:inf-sup_stable_sets}, inf-stable complete subspaces of $\pmRX$ are sets of non-expansive functions for the $\pmR$-valued weak metric given by $e_x$, so \Cref{prop-infty} can be easily adapted to this abstract setting. 
        \begin{proof}
        Case (i).  We first observe that if a net $x_\alpha\in \X$ converges to a point $x$ in the metric $\delta_s$, then $\delta(\cdot,x_\alpha)$ converges pointwise to $\delta(\cdot,x)$. To see this, note that
        \begin{equation*}
            \delta(y,x_\alpha)-\delta(y,x)
          \leq \delta(y,x)+\delta(x,x_\alpha)-\delta(y,x)= \delta(x,x_\alpha).
        \end{equation*}
            Dually, $\delta(y,x_\alpha)-\delta(y,x)\geq \delta(y,x_\alpha)-\delta(y,x_\alpha)-\delta(x_\alpha,x) = -\delta(x_\alpha,x)$. Hence, using that $(\delta(\xb,x)-\delta(\xb,x_\alpha))\in[-\delta(x,x_\alpha),\delta(x_\alpha,x)]$
          \begin{equation*}
              -\delta_s(x,x_\alpha) \le \delta(y,x_\alpha)-\delta(\xb,x_\alpha)-(\delta(y,x)-\delta(\xb,x)) \le \delta_s(x,x_\alpha).
          \end{equation*}
          Suppose now that $x_\alpha$ is an almost-geodestic and $\delta_s(\bar x, x_\alpha)$ has a bounded subnet. Then, since $(\X,\delta_s)$ has compact balls by assumption, we may assume $x_\alpha$ converges to some $x\in \X$ for $\delta_s$. Then, using the above observation, $\imath(x_\alpha)$ converges
          to $\imath(x)\in \imath(\X)$, contradicting that $h=\lim_\alpha \imath(x_\alpha)
          \in \X(\infty)= \overline{\imath(\X)}\setminus \imath(\X)$. 

        Case (ii). Assume that $\delta_s(\bar x, x_\alpha)$ has a bounded subnet, w.l.o.g.\ we restrict to this subnet. Take $\alpha \le \alpha'$, then using the triangle inequality and the almost geodesic property, for any $\epsilon>0$ and for $\alpha$ large enough
              \begin{equation*}
                  \delta_s(\bar x, x_\alpha) =\delta(x_\alpha,\xb)+\delta(\bar x, x_\alpha)\le [\delta(x_\alpha,x_{\alpha'})+\delta(x_{\alpha'},\xb)]+[\delta(\bar x, x_{\alpha'})-\delta(x_\alpha,x_{\alpha'})+\epsilon] = \delta_s(\bar x, x_{\alpha'})+\epsilon.
              \end{equation*}
              So $(\delta_s(\bar x, x_{\alpha}))_\alpha$ is almost non-decreasing, and bounded by assumption, hence it converges.%

              So, as $\delta=\frac12\delta_s$, the almost geodesic property $\delta(\bar x, x_\alpha) + \delta(x_\alpha,x_{\alpha'})\leq \epsilon+ \delta(\bar x, x_{\alpha'}) $ gives that for $\alpha\le \alpha'$ large enough, $\delta(x_\alpha,x_{\alpha'})\le 2 \epsilon$. Using again that $\delta=\frac12\delta_s$, we have that $\delta_s(x_\alpha,x_{\alpha'})\le 4 \epsilon$ for any $\alpha,\alpha'$ large enough. Since $(X,\delta_s)$ is complete, setting $\beta=\max(\alpha,\alpha')$, then 
              \begin{equation*}
                0\le \delta_s(x_\alpha,x_{\alpha'}) \le \delta_s(x_\alpha,x_{\beta})+\delta_s(x_\beta,x_{\alpha'}) \le 8\epsilon.
              \end{equation*}
              So $x_\alpha$ is Cauchy, and thus converges. We then conclude as in Case~(i).
          \end{proof}
The proof of Case~(ii) above is inspired by the argument of proof of Proposition 2.5 of~\cite{Walsh2018}. The following lemma shows that the Archimedean classes of the Busemann points cannot be maximal among the Archimedean classes of elements of the compactification $\overline{\imath(\X)}$.
        \begin{Lemma}\label{lemma-nonmaximal}
          Assume either that i) all the closed balls of $(\X,\delta_s)$ are compact or ii) $\delta=\frac12 \delta_s$ and $(\X,\delta_s)$ is complete. Then, if $h$ is a Busemann point, then there is no constant $a\in\R$ such that $h \geq a + \delta(\cdot, \bar x)$.
        \end{Lemma}
        \begin{proof}%
          By~\Cref{prop-infty}, $h$ is the limit of $\imath(x_\alpha)$ where $x_\alpha$
          is an almost-geodesic net such that $\delta_s(\bar x, x_\alpha)$ tends to infinity.
          Let us fix $\epsilon>0$, so that for all $\alpha,\beta$ large enough
          \[
          \delta(\bar x, x_\alpha) + \delta(x_\alpha,x_\beta)\leq \epsilon+ \delta(\bar x,x_\beta) .
          \]
          Suppose by contradiction that $a + \delta(\cdot, \bar x) \leq h$. Then,
          fixing $\alpha$ and $\beta$ large enough to have both the almost-geodesic inequality and $\|h-i(x_\beta)\|_\infty\le 1$, we get
          \begin{equation*}
              a-1 + \delta(x_\alpha, \bar x) \leq
          [\imath(x_\beta)](x_\alpha) = 
          \delta(x_\alpha,x_\beta)-\delta(\bar x,x_\beta)
          \leq -\delta(\bar x,x_\alpha) + \epsilon
          \end{equation*}
          and so $a-1 + \delta(x_\alpha,\bar x) + \delta(\bar x,x_\alpha) = a-1+\delta_s(\bar x ,x_\alpha )\leq \epsilon$. However, by~\Cref{prop-infty},
          the almost-geodesic $x_\alpha$ is such that $\delta_s(\bar x,x_\alpha)$ is unbounded, so we arrive
          at a contradiction.%
        \end{proof}
                \begin{Remark}
                  The compactness or completeness assumption cannot be dismissed in~\Cref{lemma-nonmaximal}. Indeed,  if $\X=\mathbb{Q}$, then, every element $z\in\R\setminus \Q$ yields an inf-irreducible horofunction $(d(\cdot,z)-d(\bar x, z))\in \Lipb(\Q)$ which is in the same Archimedean class than every map $\imath(x)$ with $x\in \X$, while having $d(\cdot,z)\notin \imath(\X)$. \Cref{lemma-nonmaximal} allows to separate the Busemann points from the distance functions $\imath(\X)$ based on the domination criterion $h \geq a + \delta(\cdot, \bar x)$. For instance for the space of 1-Lipschitz functions over a Banach space $\X$, the Busemann points are the inf-irreducible $\bk{p,\cdot}$ for $\|p\|_{\X^*}=1$, and we can thus exclusively focus on the inf-irreducible functions $e_x=\delta(\cdot,x)$ with maximal Archimedean class.
          \end{Remark}

        \begin{proof}[Proof of~\Cref{th-Lip}]\label{proof:th-lip}
          A (max,+) isomorphism $J$ of $\Lipb(\X)$ must preserve inf-irreducible elements.           Hence, by~\Cref{th-agw}, denoting by $\mathcal{B}(\X)$ the set of Busemann points
          of $\X$, $J$ must preserve the set $\R + (\imath(\X)\cup \mathcal{B}(\X))$.
          Owing to the triangle inequality, the elements of $\imath(\X)$ all belong to the same Archimedean class, which is maximal.
          By~\Cref{claim:archimd_preserved}, a (max,+) isomorphism must send a maximal Archimedean class to a maximal
          Archimedean class.  By~\Cref{lemma-nonmaximal}, the Archimedean
          class of a Busemann point is not maximal. We conclude that $J$ preserves the set $\R + \imath(\X)$.

           By \Cref{thm:max_plus_isophi_general}, we thus have \eqref{eq:charac_J}, that is $Jf=g+f\circ \phi$. We just have to characterize the properties of $g$ and $\phi$. Since $J$ sends nonexpansive maps onto nonexpansive maps, we have that, for all nonexpansive maps $f:\Y\rightarrow\R$, and all $x,x'\in\X$
        \begin{equation}\label{eq:Lip_assym1}
			Jf(x)-Jf(x')=g(x)-g(x')+f(\phi(x))-f(\phi(x'))\le \delta(x,x').
		\end{equation}
		Setting $f=\delta'(\cdot,\phi(x'))$, since $2\delta'(\phi(x),\phi(x))=\delta'(\phi(x),\phi(x))=0$, we get that for all $x,x'\in\X$
		\begin{equation}\label{eq:Lip_non_sym}
			g(x)-g(x')+\delta'(\phi(x),\phi(x'))\le  \delta(x,x').
		\end{equation}
    	Applying the same reasoning to $J^{-1}f(y)=-g(\phi^{-1}(y))+f(\phi^{-1}(y))$ for $f\in\G$ entails that for any $y,y'\in\Y$
            \begin{equation}\label{eq:Lip_non_sym2}
			-g(\phi^{-1}(y))+g(\phi^{-1}(y'))+\delta(\phi^{-1}(y),\phi^{-1}(y'))\le  \delta'(y,y').
		\end{equation}  
        Replacing $y$ by $\phi(x)$ and $y'$ by $\phi(x')$ in the latter, we obtain
        \begin{equation}\label{eq:Lip_isom_proof}
			\delta(x,x') \le g(x)-g(x')+\delta'(\phi(x),\phi(x')) \le  \delta(x,x')
		\end{equation} 
        which gives the almost isometry property. Commuting the role of $x$ and $x'$ in \eqref{eq:Lip_isom_proof} and summing, we derive that $\delta_s'(\phi(x),\phi(x'))=\delta_s(x,x')$. If $\delta$ and $\delta'$ are metrics, then $g(x)-g(x')=0$, so $g(\cdot)\equiv \alpha$ is constant.

        Conversely if \eqref{eq:max_plus_isophi_Lip} holds, then \eqref{eq:Lip_assym1} and \eqref{eq:Lip_non_sym2} hold, so $J$ maps $\Lipb(\X)$ to $\Lipb(\X)$, commutes with the (max,+) operations, and has an inverse with the same properties.
        \end{proof}

        \begin{Example}[Isomorphisms of Hilbert and Funk distances]\label{example:Funk}
          Let $C$ be a closed convex pointed cone in $\R^n$.
          Such a cone defines a partial order $\leq$ on $\R^n$ such that $x\leq y$
          iff $y-x\in C$. We define
          the (reverse) {\em Funk} ``metric'' in the interior of $C$ as follows:
          \[
          \rfunk(x,y) = \log \inf\{ \lambda>0 \mid \lambda x\geq y\} . 
          \]
It is a weak metric, see~\cite{papadopoulos,walsh}
and~\cite{funk}. The Hilbert metric on the interior of $C$ is defined by
\[
\hilbert(x,y)=\rfunk(x,y)+\rfunk(y,x).%
\]

It will be convenient to fix a linear form $p$ in the interior of the dual cone of $C$,
and to consider the cross-section $C_p = \{x\in \interior{C}\mid p(x)=1\}$.
Then,
$(x,y) \mapsto \hilbert (x,y)$ is a bona fide metric on $C_p$. Observe
that all the Hilbert balls of $C_p$ are compact, so that \Cref{th-Lip} can be applied.
The group of isometries of $(C_p,\hilbert)$
has been characterized in \cite[Theorem 1.2]{LEMMENS2011}
when $C=\R_{\geq}^n$ is the standard
orthant, so that $C_p$ is the simplex.
Let $\jmath$ denote the self-map of $(\R_{>0})^n$ such that $(\jmath(x))_i=x_i^{-1}$
for all $i\in[n]$, and $\rho$ denote the map from $\R_{>0}^n$ to
$C_p$ such that $\rho(x)=x/p(x)$. Then, every isometry of $(C_p,\hilbert)$
is either of the form $\phi(x) = \rho(DPx)$ where $D$ is a diagonal matrix
with positive diagonal entries and $P$ is a permutation matrix,
or of the form $\phi(x)=\rho(\jmath(DPx))$. Then, \Cref{th-Lip} shows
that every (max,+)-isomorphism of the space of nonexpansive maps
defined on $(C_p,\hilbert)$ is of the form~\eqref{eq:max_plus_isophi_Lip}
where the function $g$ is constant. \Cref{th-Lip} also allows
to characterize the set of (max,+)-isomorphisms of the space
of real nonexpansive maps defined on $(C_p,\rfunk)$. When $n\geq 3$,
considering $\delta=\delta'=\rfunk$ in~\eqref{eq:max_plus_isophi_Lip},
it is straightforward to check that $g(x)-g(x')+\delta'(\phi(x),\phi(x'))=\delta(x,x')$
holds for all $x,x'\in C_p$ iff $g$ is constant and $\phi$ is of the form
$\phi(x)=\rho(DPx)$ (the case in which $\phi=\rho\circ \jmath(DPx)$ is excluded owing to the ``unilateral''
character of the Funk metric).%

          \end{Example}

    \subsection{Continuous affine maps and the lower semicontinuous convex functions}\label{sec:examples_convex}

	We conclude by studying the order isomorphisms over the set of lower semicontinuous convex functions, extending the result of \cite{ArtsteinAvidan2009} to the generality of locally convex Hausdorff spaces. Our novel ingredient is that we identify in \Cref{fct:cvx_extreme_pts} the sup-irreducible points of the space of proper convex lower semicontinuous functions.

    \begin{Proposition}\label{fct:cvx_extreme_pts}
		The sup-irreducible points of the space of proper convex lower semicontinuous functions over a locally convex Hausdorff space $\X$ are the continuous affine maps $\bk{p,\cdot}+\lambda$ with $p\in\X^*$ and $\lambda\in\R$.
	\end{Proposition}
	\begin{proof} 
        $\Leftarrow$: Fix $\bp,\lambda\in\X^* \times \R$ and take $f_1,f_2$ convex l.s.c.\ such that $f_0(\cdot):=\bk{\bp,\cdot}+\lambda=\sup(f_1,f_2)$. Substracting by $f_0$, this means that the convex function $f=f_1-f_0$ is nonpositive. Convex functions smaller than constants are constants.\footnote{A short proof of this fact is that, for any $x_0,x_1$ such that $f(x_0)\le f(x_1)$, setting $x_\la=x_1+\la(x_1-x_0)$ for $\la\ge 0$, we have $x_1=\frac{1}{1+\la}x_\la+\frac{\la}{1+\la}x_0$. Hence, by convexity, $(1+\la)f(x_1)\le f(x_\la) + \la f(x_0)$ which reads as $\la(f(x_1)-f(x_0))\le f(x_\la)-f(x_1)$. The l.h.s\ diverges to $\infty$ if $f(x_0)\neq f(x_1)$, which is not compatible with $f$ being upper-bounded.} So $f_1-f_0$ and $f_2-f_0$ are both constants, and one of them must be null, so $f_0$ is thus sup-irreducible.
        
        $\Rightarrow$: Take $f$ to be a sup-irreducible proper convex l.s.c.\ function. Fixing $\bar x\in\Dom(f)$, we can assume w.l.o.g.\ $f(0)\in\R$ by replacing $f$ by $f(\cdot-\xb)$ which is also sup-irreducible. From \cite[Theorem 7.6]{aliprantis2006} we know that
        \begin{equation}\label{eq:fenchel_strict}
            f(x)=\sup_{p\in\X^*,\, \la \in\R} \{\bk{p,x}+\la\, | \, \forall \, y\in\X, \, \bk{p,y}+\la <f(y)\}=: \sup_{(p,\la)\in P^f}\bk{p,x}+\la.
        \end{equation}
        denoting by $P^f\subset \X^*\times \R$ the nonempty set over which the supremum is taken. Fix $x_0,x_1\in\X$, being possibly equal, and $(p_0,\la_0)\in P^f$, $(p_1,\la_1)\in P^f$. Assume that $p_0\neq p_1$, then by Hahn-Banach's theorem applied to the locally convex space $\X^*$ equipped with the weak-* topology, we can take a  weak-* closed halfspace $H^+$ containing $p_0$ but not $p_1$ and denote by $H^-$ the closure of its complement. This means that $H^-=\{p\in\X^* \, | \, \bk{p,x_H}\le \alpha_H\} $ for some $x_H\in\X\backslash \{0\}$ and $\alpha_H\in\R$. 
        
        Let us get rid of $\alpha_H$ first. Set $\tilde f=f-\bk{\bar p,\cdot}$ where $\bar p\in\X^*$ is such that $\bk{\bar p,x_H}= \alpha_H$. Then $\tilde f$ is sup-irreducible and $P^{\tilde f}=P^f -(\bar p, 0)$. So replacing $p_0$ by $p_0-\bar p$ and $p_1$ by $p_1-\bar p$, both are still separated, but this time by $H^-=\{p\in\X^* \, | \, \bk{p,x_H}\le 0\} $.
        
        Let $g(x):=\sup_{(p,\la)\in (H^+\times \R)\cap P^{\tilde f}}\bk{p,x}+\la$ and $h(x):=\sup_{(p,\la)\in (H^-\times \R)\cap P^{\tilde f}}\bk{p,x}+\la$. Then clearly $\tilde f=\sup(g,h)$ but $\tilde f\neq g$ and $\tilde f\neq h$. Indeed,  since  $p_0\notin H^-$, $\bk{p_0,x_H}>0$. Consequently we can fix $\beta \in\R_+^*$ such that $\beta\bk{p_0,x_H} \ge \tilde f(0)-\la_0$. For all $(p,\la)\in (H^-\times \R)\cap P^{\tilde f}$, we deduce from \eqref{eq:fenchel_strict} that $\tilde f(0) \ge \la$ and we have that $\bk{p_0,x_H} >0 \ge \bk{p,x_H} $, whence $\bk{p_0,\beta x_H}+\la_0 \ge \tilde f(0) \ge \bk{p,\beta x_H}+\la$ since $\bk{p,\beta x_H}\le 0$. Since $(p_0,\la_0)\in P^{\tilde f}$, taking the supremum in $(H^-\times \R)\cap P^{\tilde f}$ and fixing $y_0=\beta x_H$ shows that $h(y_0)\le \bk{p_0,y_0}+\la_0< \tilde f(y_0)$ so $h\neq \tilde f$. The same argument applies to $g$. Hence, $p_0=p_1$, and
        \begin{equation*}
            \tilde f(x)=\sup_{\la \in\R} \{\bk{p_0,x}+\la\, | \, \forall \, y\in\X, \, \bk{p_0,y}+\la <f(y)\}=\bk{p_0,x}-\sup_y[ \bk{p_0,y}-f(y)].
        \end{equation*}
        Consequently $\tilde f$ is a continuous affine function, and so is $f=\tilde f +\bk{\bar p,\cdot}$.
	\end{proof}

    The following theorem has previously been shown for $\X=\R^d$ in \cite[Corollary 8]{ArtsteinAvidan2009}. Its generalization to Banach spaces appeared in \cite[Theorem 1]{Iusem2015}. All proofs rely eventually on a key result of affine geometry: in dimension greater than 2 only affine maps preserve segments. It should be highlighted that it is  already stated in \cite[p670]{ArtsteinAvidan2009}, before their Remark 9, that the characterization of anti-involutions given in \Cref{cor:antiinvol_convex} could be obtained from that of order isomorphisms as in \Cref{thm:isophi_convex}. Artstein-Avidan and Milman tackle first anti-involutions and then the isomorphisms. As shown in \Cref{lem:isophi_extreme_pts}, the difference between the two is that isomorphisms send relatively-inf-irreducible elements (here the indicator functions) onto themselves, while anti-involutions send them onto the sup-irreducible elements (here the continuous affine functions). This explains the main discrepancies of our proof with that of \cite[Theorem 5]{ArtsteinAvidan2009} and its similarity with that of \cite[Theorem 1]{Iusem2015}. In constrast to \cite{Iusem2015}, an argument based on weak-* continuity enables us to avoid all together considering the bidual in the proof. Aside from this result, our proof method is different because we emphasize the role of the sup-irreducible elements to obtain the characterization. Concerning anti-involutions, \cite[Theorem 1.3]{Cheng2021} shows that reflexivity and a linear isomorphism between $\X$ and $\X^*$ are necessary and sufficient to have the Fenchel transform as anti-involution, as otherwise the bi-Fenchel is defined on the larger set $\X^{**}$. Using different proof techniques specific to the convex case, an even more general result to ours was obtained recently in \cite[Theorem 2.13]{leung2023} with points in convex subsets of locally convex Hausdorff space and values in ordered topological vector spaces.
    
	\begin{Theorem}\label{thm:isophi_convex}
		Let $\G$ be the space of proper convex lower semicontinuous functions over a locally convex Hausdorff topological vector space $\X$, then the order isomorphisms over $\G$ are affine, i.e.\ there exists a linear map $A:\X\to\X$ with w*-w* continuous adjoint, which has an inverse with w*-w* continuous adjoint, and some constants $c\in \X$, $b\in \X^*$, $d,\delta\in \R$ with $d>0$ such that, for all $f\in\G$, we have
        \begin{equation}\label{eq:isophi_convex}
            Jf(x)=\bk{b,x} +\delta +d \cdot f(Ax+c),
        \end{equation}
        and these $J$ are the only order isomorphisms over $\G$.
	\end{Theorem}
    As shown in the proof, the existence of a w*-w* continuous adjoint implies both that the linear map is w-w continuous and has closed graph. If $\X$ is Fréchet, $A$ is then continuous with continuous inverse.
	\begin{proof}
	    Using \Cref{lem:isophi_extreme_pts} and \Cref{fct:cvx_extreme_pts}, we have that $J$ maps continuous affine functions onto continuous affine functions since they are the sup-irreducible points in this case. Let $G:\X^*\times \R \rightarrow \X^*\times \R$ be the bijective function that maps the coefficients $(p,\lambda)$ of an affine function $f$ onto the coefficients $(p',\lambda')$ of $Jf$. Let us write $G(p,\la)=(G_1(p,\la), G_2(p,\la))$

    Take any two affine functions $f_0(\cdot):=\bk{p_0,\cdot}+\lambda_0$ and $f_1(\cdot):=\bk{p_1,\cdot}+\lambda_1$ for some $p_0\neq p_1 \in \X^*$, $\la_0,\la_1 \in \R$. Fix any $\lambda'_0\in\R$, w.l.o.g.\ assume that $\lambda'_0 > \lambda_0$. Since $J$ preserves the order, $J f_0 \le J(f_0+ \lambda'_0 - \lambda_0)$, both functions being affine, they must have the same $p$-coefficient, which thus does not depend on $\la$. In other words, $G_1(\cdot,\la)=G_1(\cdot,0)$.
    
    Set $f_\alpha(\cdot):=\bk{p_\alpha,\cdot}+\la_\alpha$ for $\alpha\in[0,1]$, $p_\alpha=\alpha p_1+(1-\alpha)p_0$ and $\la_\alpha=\alpha \la_1+(1-\alpha)\la_0$, then
    \begin{equation*}
        f_\alpha \le \max(f_0,f_1)
    \end{equation*}
    with equality at some point $x_0$. Consequently, since $J$ commutes with the max operation by \Cref{lem:order_implies_max}, we have that  $Jf_\alpha \le \max(Jf_0,Jf_1)$ with all functions being affine. Let us show that there exists a point $x'_0$ where there is equality.\\

    If $Jf_\alpha < \max(Jf_0,Jf_1)$, let us show there exists $\epsilon>0$ such that
    \begin{equation}\label{eq:absurd_Jfalpha}
        Jf_\alpha < Jf_\alpha+\epsilon \le \max(Jf_0,Jf_1).
    \end{equation}
    Set $f'_0=Jf_0-Jf_\alpha$, $f'_1=Jf_1-Jf_\alpha$, $F=\max(f'_0,f'_1) >0$. Denote by $G_f\subset \X\times \R$ the graph of a function $f:\X\to\R$. Let $\C:=G_{f'_0}\cap G_{f'_1}$. It is nonempty, as otherwise the hyperplanes are parallel, $G_1(p_0)=G_1(p_1)$, so $p_0=p_1$. Take $(x'_0,\la'_0)\in (G_F\cap G_{f'_0})\backslash \C$ and $(x'_1,\la'_1)\in (G_F\cap G_{f'_1})\backslash \C$, both sets being non-empty since $f'_0$ is not parallel to $f'_1$. We are down to a unidimensional reasoning in $P=(\{x'_0\}+\Sp(x'_1-x'_0))\times \R$. There clearly exists a unique $s\in (0,1)$ such that $x'_s=(1-s)x'_0+sx'_1$ satisfies $(x'_s,F(x'_s))\in \C$. If w.l.o.g.\ $F(x'_0)<F(x'_s)$, then since $\C$ is a vector space, $x(\cdot):\beta\in (-\infty,s]\mapsto (1-\beta)x'_0+\beta x'_1$ is such that $F(x(\beta))=f'_0(x(\beta))$. Taking the limit $\beta\to-\infty$, by linearity of $f'_0$, we obtain $F(x(\beta))\to -\infty$, while $F$ is lower bounded. Hence we must have $F(x_s)\le \min(F(x'_0),F(x'_1))$, and to estimate $\inf_\X F$, we can focus exclusively on $\C$. Take now $(x'_0,\la'_0)\in \C$ and $(x'_1,\la'_1)\in \C$ with $x'_0\neq x'_1$, and define $P$ as before. If w.l.o.g.\ $\la'_0>\la'_1$, then setting $x_\beta:=(1-\beta)x'_0+\beta x'_1$, since $F_{|P}=f'_0=f'_1$ is affine, we obtain $\lim_{\beta\to+\infty}F(x_\beta)=-\infty$. So $\la'_0=\la'_1$, i.e.\ $\C\subset \X\times\{\epsilon\} $ for some $\epsilon\in\R$. Since $0<F$, we necessarily have $\epsilon>0$.\\    
    
    Having shown that $Jf_\alpha < \max(Jf_0,Jf_1)$ implies \eqref{eq:absurd_Jfalpha}, as $G$ is bijective, $J^{-1}$ commutes with the max, and as the $p$-coefficient of $J^{-1}(Jf_\alpha+\epsilon)$ does not depend on $\epsilon$, we have that $J^{-1}(Jf_\alpha+\epsilon)=\bk{p_\alpha,\cdot}+\la_{\alpha,\epsilon}$ for some $\la_{\alpha,\epsilon}$. We have furthermore that $\la_{\alpha,\epsilon}>\la_\alpha$ since $G_2(\_,\cdot)$ is bijective and increasing in $\la$. Hence, applying $J^{-1}$ to \eqref{eq:absurd_Jfalpha}, we derive that
    \begin{equation*}
        f_\alpha=\bk{p_\alpha,\cdot}+\la_\alpha < \bk{p_\alpha,\cdot}+\la_{\alpha,\epsilon}\le \max(f_0,f_1)
    \end{equation*}
    which contradicts the existence of $x_0$.
    
    Applying a convex relaxation to the r.h.s.\ of our inequality $Jf_\alpha \le \max(Jf_0,Jf_1)$, with equality at one point, we can consequently write that
    \begin{equation*}
        0= \min_{x\in\X} \max_{s,t\ge 0, \, s+t=1} s(Jf_0(x)-Jf_\alpha(x))+t(Jf_1(x)-Jf_\alpha(x)).
    \end{equation*}
    Using Sion's minimax theorem \cite[Theorem 4.2']{Sion1958} over these continuous affine functionals, we have strong duality, so 
    \begin{equation*}
        0=  \max_{s,t\ge 0, \, s+t=1} \min_{x\in\X} s(Jf_0(x)-Jf_\alpha(x))+t(Jf_1(x)-Jf_\alpha(x))
    \end{equation*}
    Let $s,t$ be maximizing in the above. The function minimized is affine, the minimization over $x$ implies that it is constant. It is thus null since the value of the max is zero. In other words, the two coefficients determining the affine function $Jf_\alpha$ are jointly a convex combination of those of $Jf_0$ and $Jf_1$.
    
    We have just shown that $G$ maps intervals onto intervals. Arguing as in \cite[Remark 6, p668]{ArtsteinAvidan2009} by applying a result of affine geometry, extended in \cite[Corollary 3]{Iusem2015} to arbitrary vector spaces, since $\X\times\R$ has dimension larger than two, $G$ is thus affine and bijective, i.e.\ we can fix a linear invertible $A:\X^*\to\X^*$, a linear form $C:\X^*\to\R$, and constants $g,b\in \X^*$, $d,\delta\in \R$ such that
    \begin{equation*}
        J(\bk{p,\cdot}+\lambda)=\bk{Ap+g\lambda+b,\cdot}+Cp+d\lambda+\delta.
    \end{equation*}
    We have already shown before that $g=0$ as there was no dependence in $\la$ in the $p$-term. Applying $J$ to the constant functions $0\le 1$, we obtain that 
    \begin{equation*}
        \bk{b,\cdot}+\delta\le \bk{b,\cdot}+d+\delta.
    \end{equation*}
    Hence $d\ge 0$. Since $G$ is bijective in $\la$, we have that $d>0$ necessarily. It remains to study the continuity of the maps $A$ and $C$. Since $\G$ contains the $\{\delta_x^\top+\la\}_{x\in\X,\la\in\pmR}$, we derive from the intermediary result \eqref{eq:isophi_lsc_inter} in the proof of \Cref{thm:max_isophi_lsc} that there exists a lower semicontinuous $g:\X\times \R\to\R$, with, for all $x\in\X$, $g(x,\cdot):\R\to\R$ continuous (being increasing and bijective), and a bijective $\phi:\X\to\X$ such that
    \begin{equation*}
        J f(x)= g(x,f(\phi(x))),\, \forall x\in\X, \, f\in\G
    \end{equation*}
    In particular for $f=\bk{p,\cdot}+\lambda$, we obtain that
    \begin{equation*}
        \bk{Ap+b,x}+Cp+d\lambda+\delta= g(x,\bk{p,\phi(x)}+\lambda),\, \forall x\in\X, \, p\in\X^*.
    \end{equation*}
    Take $x=0$ and $\la=0$, then $Cp+\delta= g(0,\bk{p,\phi(0)})$ so $C$ is weak-* continuous as $g(0,\cdot)$ is continuous, i.e.\ there exists $c \in\X$ such that $Cp=\bk{p,c}$. Moreover
    \begin{equation*}
        \bk{Ap+b,x}=-g(0,\bk{p,\phi(0)})+ g(x,\bk{p,\phi(x)}),\, \forall x\in\X, \, p\in\X^*,
    \end{equation*}
    so, considering nets, we have that $A:\X^*\to\X^*$ is weak-* continuous as $g(x,\cdot)$ is continuous. Hence $A$ has an adjoint $A^*:\X\to\X$ with closed graph. Indeed, by weak-* continuity, for all $x\in\X$, there exists $y_x\in\X$ such that for all $p\in\X^*$, $\bk{Ap,x}=\bk{p,y_x}$. Set $A^* (x):=y_x$, since $\bk{Ap,x}=\bk{p,A^* (x)}$ and $\X^*$ separates the points of the locally convex $\X$, $A^*$ is linear. Take a net $(x_\alpha,y_\alpha)_\alpha$ with $y_\alpha=A^*x_\alpha$ converging to some $(x,y)\in\X\times \X$, then $\bk{p,y}=\bk{Ap,x}=\bk{p,A^*x}$, so $y=A^*x$ and the graph of $A^*$ is closed. Similarly if $(x_\alpha)_\alpha$ converges weakly to some $x\in\X$, then $(A^*x_\alpha)_\alpha$ converges weakly to $A^*x$.
    
    Fix any convex l.s.c.\ $f$. Then
    \begin{align}
        Jf(x) &\stackrel{\eqref{eq:fenchel_strict}}{=} \sup_{(p,\la)\in P^f} J(\bk{p,\cdot}+\la)(x) \nonumber \\
        &= \sup_{(p,\la)\in P^f} \bk{Ap+b,x}+\bk{p,c}+d\la+\delta \nonumber\\
        &= \bk{b,x} +\delta +d\cdot\sup_{(p,\la)\in P^f} [\bk{p,\frac1d(A^*x+c)}+\la]\stackrel{\eqref{eq:fenchel_strict}}{=}\bk{b,x} +\delta +d \cdot f(\frac1d(A^*x+c)).\label{eq:Jf_sci}
    \end{align}
    Renaming the parameters, we have obtained the desired result \eqref{eq:isophi_convex}.

    Conversely, take any $J$ satisfying \eqref{eq:isophi_convex}. Such $J$ is clearly order preserving and, reversing the reasoning which led to \eqref{eq:Jf_sci} and using that $A$ has a weak-* continuous adjoint, $Jf$ is the supremum of lower semicontinuous convex functions, hence $Jf\in\G$. If we set $J^1f(x)=-\frac1d\bk{b,A^{-1}(x-c)} -\frac{\delta}{d} +\frac1d \cdot f(A^{-1}(x-c))$, it satisfies $J^1Jf=f$.

	\end{proof}
    Applying \Cref{prop:isophi_anti-invol}, we know that, for $\G=\Rg(B)$ for some $\bB$ as in \eqref{eq:RgB_tpsd}, the anti-involutions are of the form $\bB J$ with order isomorphisms $J:\G\to \G$ such that $\bB J \bB J=\Id_{\G}$. We can equivalently consider $J \bB$. In the present case, $\bB$ is the Fenchel transform, and we thus recover \cite[Theorem 5]{ArtsteinAvidan2009} for $\X=\R^d$ as a special case and derive a straightforward extension to semi-convex functions, i.e.\ the functions $f$ such that $f+\frac{\alpha}{2}\|\cdot\|^2$ is convex l.s.c.
    \begin{Corollary}\label{cor:antiinvol_convex}
	   Let $\G$ be the space of proper convex lower semicontinuous functions over a reflexive Banach space $\X$, assumed to be linearly isomorphic to its dual $\X^*$. Then the anti-involutions over $\G$ are of the form
        \begin{equation}\label{eq:anti_invol_convex}
            Tf(x)=\bk{Kc,x} +\delta + f^* (K(Ax+c)).
        \end{equation}
        with $A\in \Lcal(\X)$ invertible, $c\in \X$, $\delta\in\R$ and $K^{-1}A^{-\top}KA=\Id_\X$ where $K:\X \rightarrow \X^*$ is the duality operator and $A^{-\top}$ is the inverse of the adjoint of $A$.%

        Let $\alpha\in\R$. If $\G$ is instead the space of proper $\alpha$-semiconvex lower semicontinuous functions, then the anti-involutions over $\G$ are of the form  \begin{equation}\label{eq:anti_invol_semiconvex}
            Tf(x)=-\frac{\alpha}{2}\|x\|_\X^2+\bk{Kc,x} +\delta +(f+\frac{\alpha}{2}\|\cdot\|_\X^2)^* (K(Ax+c)),
        \end{equation}
        with $A,K,c$ as above. Such $J$ have an alternative expression based on the inf-convolution, using that $(f+\frac{\alpha}{2}\|\cdot\|_\X^2)^*=f^*\square(\frac{1}{2\alpha}\|\cdot\|_{\X^*}^2)$.
	\end{Corollary}
 \begin{proof}
     For convex functions, applying \Cref{prop:isophi_anti-invol}, since the Fenchel transform is an anti-involution with these assumptions on $\X$, see e.g.\ \cite[Theorem 1.3]{Cheng2021}, we have that $Tf(x)=\bk{b,x}+\delta +d \cdot f^* (K(Ax+c))$. The operator $TT$ is a (max,+)-isomorphism since $TT=\Id_\G$. However $TT(f+\la)=TTf+d^2\la$, so necessarily $d=1$. Moreover $\{\bk{p,\cdot}+\la\}_{p\in\X^*,\la\in\R}$ sup-generates $\G$, so $TT=\Id_\G$ iff, for all $p\in\X^*$, $TT\bk{p,\cdot}=\bk{p,\cdot}$. Fix $p\in\X^*$, we have that $T\bk{p,\cdot}(x)=\bk{b,x}+\delta+\delta^\top_{A^{-1}(K^{-1}p-c)}$, and $TT\bk{p,\cdot}(x)=\bk{b,x}+ \bk{K(Ax+c)-b,A^{-1}(K^{-1}p-c)}$. Identifying the ``monomials'' in $p$ and $x$, we get that $b=Kc$ and $K^{-1}A^{-\top}KA=\Id_\X$ are necessary and sufficient conditions.

     For the semi-convex functions, we just use \Cref{lem:extension_order} since by definition a function $f$ is $\alpha$-semiconvex iff $f+\frac{\alpha}{2}\|\cdot\|_\X^2$ is convex.
 \end{proof}

    \paragraph{Acknowledgments.} The authors deeply thank Sebastián Tapia-García for the fruitful discussions and his many insights on topology; as well as a very thorough anonymous reviewer for their detailed reading which improved the paper. Most of the article was written while Pierre-Cyril Aubin-Frankowski was hosted at TU Wien, being funded by the FWF project P 36344-N.

        \section*{Appendix}
        \subsection*{Proof of \Cref{thm:sym_anti-involution}}\label{proof:sym_anti}
        \ref{it_rgB}$\Rightarrow$\ref{it_antiinvol}. We use the classical property that $\bB^\circ \bB \bB^\circ = \bB^\circ$, see e.g.\  \cite[p.3]{Akian04setcoverings}, so $\bB$ is an anti-isomorphism from $\Rg(B^\circ)$ onto $\Rg(B)$.

		\noindent	\ref{it_antiinvol}$\Rightarrow$\ref{it_rgB}. By \Cref{lem:order_implies_max}, anti-isomorphisms $\bF$ send relative infima $\inf\nolimits^\F$ onto suprema. Define $\tilde{F}:\pmRY\rightarrow \pmRX$ by, for any $f\in\pmRY$, $\tilde{F}(f):=\bF(\inf\nolimits^\F f) \in \G$. Similarly, set $\tilde{F}^1:\pmRX\rightarrow \pmRY$ as for any $g\in\pmRX$, $\tilde{F}^1(g):=\bF^1(\inf\nolimits^\G g) \in \F$.  
  
    The function $\tilde{F}$ is automatically an extension of $\bF$ and $\tilde{F}^1\bF=\Id_\F$. Moreover, since $\bF$ is bijective and order reversing, $\tilde{F}(x\mapsto \pm \infty)=\bF(x\mapsto \pm \infty)=\mp \infty$. Fix $f\in\pmRY$ and $\lambda\in\R$, then, as $\bF$ commutes with the addition of scalars, $\tilde{F}(f+ \lambda)=\bF(\lambda+\inf\nolimits^\G f)=\tilde{F}(f)-\lambda$. Let $(f_\alpha)_{\alpha\in \Asc}\in(\pmRY)^{\Asc}$. By definition, $\inf\nolimits^\F(\inf_{\alpha\in \Asc} f_\alpha)=\inf\nolimits^\F_\alpha (\inf\nolimits^\F f_\alpha)$. Hence
		\begin{equation*}
			\tilde{F}(\inf_{\alpha\in \Asc} f_\alpha) = \bF(\inf\nolimits^\F_\alpha (\inf\nolimits^\F f_\alpha))=\sup_\alpha \bF(\inf\nolimits^\F f_\alpha)= \sup_\alpha\tilde{F}(f_\alpha). 
		\end{equation*}
		Hence $\tilde{F}$ is continuous in the sense of \Cref{def:rmax-sesqui-lin} and commutes with the addition of $\pmR$-scalars. We can now apply \Cref{thm:singer} to derive a kernel $b$ associated with $\tilde{F}$. Since $\tilde{F}^1\bF=\Id_\F$, we have that $\F= \tilde{F}^1\Rg(\tilde{F})$, so $\G=\Rg(\tilde{F})$. As $\bF$ is an anti-isomorphism, 
		\begin{displaymath}
			\tilde{F}^1\tilde{F}(f)=\bF^1\bF(\inf\nolimits^\F f)=\inf\nolimits^\F f\le f.
		\end{displaymath}
		We have thus recovered one of the characterizations of dual Galois connections over $\pmRY$, as given in \cite[p.3, Eq.(2a)]{Akian04setcoverings}. \cite[Theorem 2.1, Example 2.8]{Akian04setcoverings} then allows to conclude, since the dual connection is simply the transpose $\tilde{F}^\top$ \cite[Theorem 8.4]{singer1997}.

        If $\X=\Y$ and $\bF$ is an anti-involution, then $\tilde{F}$ is equal to its dual connection $\tilde{F}^1=\tilde{F}^\top$, so the kernel $b$ is symmetric.

    \bibliographystyle{alphaabbr}
	\bibliography{characterization_anti-invol}	

\end{document}